\documentclass[12pt]{article}

\usepackage{mathptmx,amsmath,amssymb,amsfonts,amsthm}
\usepackage{graphicx,color,xcolor,pictex,epstopdf,url,algorithm,algorithmic,caption,endnotes}
\usepackage{comment}

\usepackage{etoolbox} 
\usepackage{wrapfig} %
\makeatletter
\patchcmd{\@makechapterhead}{\large}{\normalsize}{}{}
\patchcmd{\@makeschapterhead}{\normalsize}{\normalsize}{}{}
\makeatother
\usepackage[textwidth=6.0in,textheight=8.9in,left=0.75in,right=0.75in,top=0.75in,bottom=0.75in]{geometry}
\usepackage[section]{placeins}
\makeatletter
\g@addto@macro\normalsize{\setlength\abovedisplayskip{4pt}}
\g@addto@macro\normalsize{\setlength\belowdisplayskip{4pt}}
\makeatother
\newtheorem{theorem}{Theorem}[section]
\newtheorem{lemma}{Lemma}[section]
\newtheorem{corollary}{Corollary}[section]
\newtheorem{proposition}{Proposition}[section]

\usepackage{parskip}
\let\oldref\ref
\renewcommand{\ref}[1]{(\oldref{#1})}  
\renewcommand{\eqref}[1]{(\oldref{#1})} 

\newbox\boxaddrone \newbox\boxaddrtwo

\input colordvi
\font\tenrm=cmr10
\font\teni=cmmi10 \skewchar\teni='177
\font\tensy=cmsy10 \skewchar\tensy='60
\font\tenex=cmex10
\font\tenit=cmti10
\font\tensl=cmsl10
\font\tenbf=cmbx10
\font\tentt=cmtt10
\font\ninerm=cmr9
\font\ninei=cmmi9 \skewchar\ninei='177
\font\ninesy=cmsy9 \skewchar\ninesy='60
\font\nineit=cmti9
\font\ninesl=cmsl9
\font\ninebf=cmbx9
\font\ninett=cmtt9
\font\eightrm=cmr8
\font\eighti=cmmi8 \skewchar\eighti='177
\font\eightsy=cmsy8 \skewchar\eightsy='60
\font\eightit=cmti8
\font\eightsl=cmsl8
\font\eightbf=cmbx8
\font\eighttt=cmtt8
\font\sevenrm=cmr7
\font\seveni=cmmi7 \skewchar\seveni='177
\font\sevensy=cmsy7 \skewchar\sevensy='60
\font\sevenbf=cmbx7
\font\sevenit=cmmi7
\font\sevensl=cmmi7
\font\seventt=cmr7
\font\sixrm=cmr6
\font\sixi=cmmi6 \skewchar\sixi='177
\font\sixsy=cmsy6 \skewchar\sixsy='60
\font\sixbf=cmbx6
\font\fiverm=cmr5
\font\fivei=cmmi5 \skewchar\fivei='177
\font\fivesy=cmsy5 \skewchar\fivesy='60
\font\fivebf=cmbx5
\def\tenpoint{\def\rm{\fam0\tenrm}%
        \textfont0=\tenrm \scriptfont0=\sevenrm \scriptscriptfont0=\fiverm
        \textfont1=\teni \scriptfont1=\seveni \scriptscriptfont1=\fivei
        \textfont2=\tensy \scriptfont2=\sevensy \scriptscriptfont2=\fivesy
        \textfont3=\tenex \scriptfont3=\tenex \scriptscriptfont3=\tenex
        \def\it{\fam\itfam\tenit}%
        \textfont\itfam=\tenit
        \def\sl{\fam\slfam\tensl}%
        \textfont\slfam=\tensl
        \def\bf{\fam\bffam\tenbf}%
        \textfont\bffam=\tenbf \scriptfont\bffam=\sevenbf
                \scriptscriptfont\bffam=\fivebf
        \def\tt{\fam\ttfam\tentt}%
        \textfont\ttfam=\tentt
        \normalbaselineskip=12pt%
        \let\sc=\eightrm        
        \setbox\strutbox=\hbox{\vrule height8.5pt depth3.5pt width0pt}%
        \normalbaselines\rm}
\def\ninepoint{\def\rm{\fam0\ninerm}%
        \textfont0=\ninerm \scriptfont0=\sixrm \scriptscriptfont0=\fiverm
        \textfont1=\ninei \scriptfont1=\sixi \scriptscriptfont1=\fivei
        \textfont2=\ninesy \scriptfont2=\sixsy \scriptscriptfont2=\fivesy
        \textfont3=\tenex \scriptfont3=\tenex \scriptscriptfont3=\tenex
        \def\it{\fam\itfam\nineit}%
        \textfont\itfam=\nineit
        \def\sl{\fam\slfam\ninesl}%
        \textfont\slfam=\ninesl
        \def\bf{\fam\bffam\ninebf}%
        \textfont\bffam=\ninebf \scriptfont\bffam=\sixbf
                \scriptscriptfont\bffam=\fivebf
        \def\tt{\fam\ttfam\ninett}%
        \textfont\ttfam=\ninett
        \normalbaselineskip=11pt%
        \let\sc=\sevenrm        
        \setbox\strutbox=\hbox{\vrule height8pt depth3pt width0pt}%
        \normalbaselines\rm}
\def\eightpoint{\def\rm{\fam0\eightrm}%
        \textfont0=\eightrm \scriptfont0=\sixrm \scriptscriptfont0=\fiverm
        \textfont1=\eighti \scriptfont1=\sixi \scriptscriptfont1=\fivei
        \textfont2=\eightsy \scriptfont2=\sixsy \scriptscriptfont2=\fivesy
        \textfont3=\tenex \scriptfont3=\tenex \scriptscriptfont3=\tenex
        \def\it{\fam\itfam\eightit}%
        \textfont\itfam=\eightit
        \def\sl{\fam\slfam\eightsl}%
        \textfont\slfam=\eightsl
        \def\bf{\fam\bffam\eightbf}%
        \textfont\bffam=\eightbf \scriptfont\bffam=\sixbf
                \scriptscriptfont\bffam=\fivebf
        \def\tt{\fam\ttfam\eighttt}%
        \textfont\ttfam=\eighttt
        \normalbaselineskip=9pt%
        \let\sc=\sixrm  
        \setbox\strutbox=\hbox{\vrule height7pt depth2pt width0pt}%
        \normalbaselines\rm}
\def\sevenpoint{\def\rm{\fam0\sevenrm}%
        \textfont0=\sevenrm \scriptfont0=\fiverm \scriptscriptfont0=\fiverm
        \textfont1=\seveni \scriptfont1=\fivei \scriptscriptfont1=\fivei
        \textfont2=\sevensy \scriptfont2=\fivesy \scriptscriptfont2=\fivesy
        \textfont3=\tenex \scriptfont3=\tenex \scriptscriptfont3=\tenex
        \def\it{\fam\itfam\sevenit}%
        \textfont\itfam=\sevenit
        \def\sl{\fam\slfam\sevensl}%
        \textfont\slfam=\sevensl
        \def\bf{\fam\bffam\sevenbf}%
        \textfont\bffam=\sevenbf \scriptfont\bffam=\fivebf
                \scriptscriptfont\bffam=\fivebf
        \def\tt{\fam\ttfam\seventt}%
        \textfont\ttfam=\seventt
        \normalbaselineskip=8pt%
        \let\sc=\fiverm  
        \setbox\strutbox=\hbox{\vrule height6pt depth2pt width0pt}%
        \normalbaselines\rm}

\input pictex
\newbox\figurelegendone
\newbox\figurelegendtwo
\newbox\figurelegendthree
\newbox\figureone
\newbox\figuretwo
\newbox\figurethree
\newbox\figurefour
\newbox\figurefive
\newbox\figuresix
\newbox\figureseven
\newbox\figureeight
\newbox\figurenine
\newbox\figureten
\newbox\figureeleven
\newbox\figuretwelve
\newbox\figurethirteen
\newbox\figurefourteen
\newbox\figurefifteen
\newbox\textone
\newdimen\xfiglen \newdimen\yfiglen  \newdimen\textlen
\definecolor{darkblue}{RGB}{16,30,144} 

\begin{document}

\title{ On an inverse problem of nonlinear imaging with fractional damping}
\author{Barbara Kaltenbacher\footnote{
Department of Mathematics,
Alpen-Adria-Universit\"at Klagenfurt.
barbara.kaltenbacher@aau.at.}
\and
William Rundell\footnote{
Department of Mathematics,
Texas A\&M University,
Texas 77843. 
rundell@math.tamu.edu}
}
\date{\vskip-3ex}
 \maketitle  

\begin{abstract}
This paper considers the attenuated Westervelt equation in pressure formulation.
The attenuation is by various models proposed in the literature and
characterised by the inclusion of non-local operators that give power law
damping as opposed to the exponential of classical models.
The goal is the inverse problem of recovering a spatially dependent coefficient
in the equation, the parameter of nonlinearity $\kappa(x)$,
in what becomes a nonlinear hyperbolic equation with nonlocal terms.
The overposed measured data is a time trace taken on a subset of the domain
or its boundary.
We shall show injectivity of the linearised map from $\kappa$ to the overposed
data used to recover it and from this basis develop and analyse Newton-type
schemes for its effective recovery.

\end{abstract}

\leftline{\small \qquad\qquad
{\bf Keywords:} 
Inverse problems, damped nonlinear wave equation, ultrasound.}
\smallskip
\leftline{\small \qquad\qquad 
\textbf{{\textsc ams}} {\bf classification:} 35R30, 35R11, 35K55, 78A46}

\section{Introduction}\label{sec:introduction}

The problem of nonlinearity parameter imaging with ultrasound in lossy media amounts to identification of the space dependent coefficient $\kappa(x)$ for the attenuated Westervelt equation
in pressure formulation
\begin{equation}\label{eqn:Westervelt_init_D_intro}
\begin{aligned}
&\bigl(1-\kappa(x)v^2\bigr)_{tt}-c_0^2\triangle v + D v = r \mbox{ in }\Omega\times(0,T)\\
v&=0 \mbox{ on }\partial\Omega\times(0,T)\;;\qquad
v(0)=0, \quad v_t(0)&=0 \;\mbox{ in }\;\Omega
\end{aligned}
\end{equation}
from observations 
Here $c_0>0$ is the wave speed (possibly space dependent as well), 
and $D v$ a damping term that will be specified below.
For simplicity we impose homogeneous Dirichlet boundary conditions here, 
but the ideas in this paper extend to more realistic boundary conditions,
such as absorbing boundary conditions for avoiding spurious reflections
and/or inhomogeneous Neumann boundary conditions for modelling excitation via,
e.g., some transducer array.
Note that the excitation here is modelled by an interior source $r$,
and we refer to a discussion on this in \cite{KaltenbacherRundell:2021c}.

By letting our $v$ equation satisfy  boundary
(or possibly interior) observations we obtain an inverse problem for the
recovery of $\kappa$.
These measurements will be taken to be
\begin{equation}\label{eqn:obs_intro}
g(t)=v(x_0,t) \quad \mbox{ or } \quad g(x,t)=v(x,t) \,,\ x\in\Sigma\,,t\in(0,T).
\end{equation}
at some point $x_0$ or on some surface $\Sigma$ contained in $\overline{\Omega}$.

The inverse problem represented by
equations~\eqref{eqn:Westervelt_init_D_intro} and \eqref{eqn:obs_intro}
is challenging on at least three counts.
First, the underlying model equation is nonlinear and in fact
the nonlinearity occurs in the highest order term.
Second, the unknown coefficient $\kappa(x)$ is directly coupled to this term
and third, is spatially varying whereas the data $g(t)$ is in the
``orthogonal'' time direction and this is well known to lead to
several ill-conditioning of the inversion of the map from data to unknown.

The use of ultrasound is a well-established protocol in the imaging of
human tissue and, besides the classical sonography methodology,
there exist several novel imaging principle, such as harmonic imaging
or nonlinearity imaging. 
The latter
\cite{Bjorno1986, BurovGurinovichRudenkoTagunov1994, Cain1986, 
IchidaSatoLinzer1983, VarrayBassetTortoliCachard2011, ZhangChenGong2001, ZhangChenYe1996}, 
relies on tissue-dependence, hence spatial variation of the parameter of
nonlinearity $B/A$ that is contained in $\kappa$.
See, e.g.,
It thus inherently needs a nonlinear acoustic model as an underlying
{\sc PDE} and we refer to e.g., the review
\cite{Kaltenbacher:reviewNonlinearAcoustics:2015}
and the references therein for a brief derivation of the fundamental
nonlinear acoustic equations.
The quantity of interest from an imaging perspective
is the coefficient $\kappa$ and its recovery in the case when the damping term
was $D u = b\triangle u_t$ was the subject of \cite{KaltenbacherRundell:2021c}.
This is the classical formulation of damping being proportional to velocity
but there are may alternative models that are prominent in the literature.
We mention some of these in the next section but the main change is
the incorporation of non-local terms involving either fractional
derivatives in time or modifying the operator $(-\triangle)$ to 
have the Laplacian raised to a fractional power $(-\triangle)^\beta$.
These have the effect of ameliorating the exponential decay of
the solution, by a fractional exponent in the frequencies in the case
of $(-\triangle)^\beta$ and by a power law decay in the case of a fractional
time derivative.
The use of such operators in inverse problems is now well documented in
the literature (see, \cite{JinRundell:2015} and in particular, for
the wave equation in \cite{KaltenbacherRundell:2021b}).

In this paper we will provide analysis for the forwards problem and in
particular regularity and well-posedness for the coefficient-to-state map
$G:\kappa\mapsto v$ where $v$ solves \eqref{eqn:Westervelt_init_D_intro}.
The forwards map is defined by $F(\kappa)=\mbox{tr}_\Sigma v$,
where $\mbox{tr}_\Sigma v$ denotes the time trace of the space--and--time
dependent function $v:(0,T)\times \Omega$ at the observation surface
$\Sigma$ (which may also just be a single point $\Sigma=\{x_0\}$).
Its linearization at $\kappa=0$ is $F'(0)$ and we will prove an injectivity
result which will both show local uniqueness and pave the way for 
the use of Newton's method which we will analyse in section~\oldref{sect:inject}
then show  reconstructions based on this approach in
section~\oldref{sec:reconstructions}.

\section{The Imaging Problem}\label{sec:imaging}

As already mentioned in the introduction, the inverse problem under
consideration is to recover the space dependent coefficient $\kappa(x)$ in
the attenuated Westervelt equation which can also be written in the form
\begin{equation}\label{eqn:Westervelt_init_D}
\begin{aligned}
u_{tt}+c^2\mathcal{A} u + D u &= \kappa(x)(u^2)_{tt} + r \mbox{ in }\Omega\times(0,T)\\
u(0)=0, \quad u_t(0)&=0 \mbox{ in }\Omega
\end{aligned}
\end{equation}
from observations 
$g(x,t)\!=\!u(x,t)$, $x\in\Sigma,\; t\in\!(0,T)$ where 
$\Sigma$ may be a single point, typically located on $\partial\Omega$.
Here $c>0$ is the constant mean wave speed, and $\mathcal{A}=-(c_0(x)^2\!/c^2) \triangle $ contains the possibly spatially varying coefficient $c_0(x)>0$ and is equipped -- for simplicity -- with homogeneous boundary conditions.
Moreover $r=r(x,t)$ is a known source term modelling excitation of the acoustic wave by a transducer array, see \cite{KaltenbacherRundell:2021b}.
 
%
The damping operator $D$ appearing in \eqref{eqn:Westervelt_init_D} is a differential operator containing space and/or time derivatives.
Classically, $D$ will consist of integer derivatives, typical examples being $D=\mathcal{A} \partial_t$ or $D=\partial_t$ often referred to as strong and weak damping, respectively.
We here list some of the (due to experimentally found power law frequency dependence) practically relevant fractional damping models, that we have already discussed in \cite{KaltenbacherRundell:2021b} in a different imaging context, namely for the inverse {\sc pat/tat} problem:
\goodbreak

\smallskip
\noindent
{\bf Time fractional models:}

\begin{itemize}
\item 
Caputo-Wismer model \cite{Caputo:1967}, \cite[eq. (5)]{Wismer:2006},
called Kelvin wave equation in  \cite[eq. (19)]{CaiChenFangHolm_survey2018}
\begin{equation}\label{eqn:D_timefrac}
D = b \mathcal{A} \partial_t^\alpha  
\end{equation}
where typically $\alpha\in[0,1]$. 
\item 
(Modified) Szabo model  \cite{Szabo:1994}, \cite[eq. (42)]{CaiChenFangHolm_survey2018}
\begin{equation}\label{eqn:D_timefracS}
D = b \partial_t^{\alpha+2}  
\end{equation}
where $\alpha\in[-1,1]$, $b\geq0$. 
\item 
Fractional Zener (combined Caputo-Wismer-Szabo) model \cite{HolmNaesholm:2011,Mainardi:2010}, \cite[eq. (30)]{CaiChenFangHolm_survey2018}
\begin{equation}\label{eqn:D_timefracZ}
D = b_1 \mathcal{A} \partial_t^{\alpha_1}  + b_2 \partial_t^{\alpha_2+2} 
\end{equation}
where $\alpha_1\geq\alpha_2\in[0,1]$, $b_1\geq b_2 c^2$, cf. \cite[Section III.B]{HolmNaesholm:2011}. 
\end{itemize}
In these models $\partial_t^\alpha$ denotes the Djrbashian-Caputo fractional time derivative, which here, due to the homogeneous initial conditions, coincides with the Riemann-Liouville one.

\smallskip
\noindent
{\bf Space fractional models:}
\begin{itemize}
\item 
Chen-Holm model \cite[eq. (21)]{ChenHolm:2004} 
\begin{equation}\label{eqn:D_spacefracCH}
D = b \mathcal{A}^\beta 
\end{equation}
typically with $\beta\in[0,1]$, where Kelvin-Voigt damping is recovered when
$\beta=1$. 
\item Treeby-Cox model \cite[eq. (28)]{TreebyCox:2010} 
\begin{equation}\label{eqn:D_spacefracTC}
D = b_1 \mathcal{A}^\beta\partial_t   + b_2 \mathcal{A}^{\beta+1/2} 
\end{equation}
typically $\beta\in[0,1]$, which is an extension of the former.
\end{itemize}
Here we use the spectral definition of the Laplacian which coincides with the Riesz version on $\mathbb{R}^d$; however, they differ in case of bounded $\Omega$.

\smallskip

In this paper, we will focus on two damping models namely (a) a combination
of \eqref{eqn:D_timefrac} and \eqref{eqn:D_spacefracCH}, since we find it
interesting to investigate the interplay of space- and time-fractional
derivatives and its influence on the ill-posedness of the inverse problem;
(b) \eqref{eqn:D_timefracZ} as it contains higher than second order time derivatives which are in case $\alpha_2=1$ known to make the equation behave wave-like (finite speed of propagation) in spite of the damping, which is expected to influence the degree of ill-posedness of the inverse problem as well.
Thus we here focus on the two damping models
$$
\begin{aligned}
D &= b \mathcal{A}^\beta \partial_t^\alpha 
\quad \mbox{ (combination of Caputo-Wismer and Chen-Holm model -- CWCH)}
\\
D &= b_1 \mathcal{A} \partial_t^{\alpha_1}  + b_2 \partial_t^{\alpha_2+2} 
\quad \mbox{ (fractional Zener -- FZ)}\,.
\end{aligned}
$$

\section{Analysis of the forward problem}

In this section, we consider well-posedness of the initial value problem for the parameter-to-state map $G:\kappa\mapsto u$ where $u$ solves \eqref{eqn:Westervelt_init_D}
and its linearisation 
$z=G'(\kappa)\underline{d\kappa}$  
\begin{equation}\label{eqn:Westervelt_init_D_lin}
\begin{aligned}
(1-2\kappa u)z_{tt}+c^2\mathcal{A} z + D z - 4 \kappa u_t\, z_t - 2\kappa u_{tt} \, z
&= 2\underline{\delta\kappa}(u \,u_{tt} + u_t^2) \mbox{ in }\Omega\times(0,T)\\
z(0)=0, \quad z_t(0)&=0 \mbox{ in }\Omega
\end{aligned}
\end{equation}
for given $\kappa$ and $\underline{\delta\kappa}$, respectively, in the context of several damping models $D$.
In order to prove Fr\'{e}chet differentiability, we will also have to consider the difference $v=G(\tilde{\kappa})-G(\kappa)=\tilde{u}-u$, which solves
\begin{equation}\label{eqn:v}
\begin{aligned}
(1-2\kappa u)v_{tt}&+c^2\mathcal{A} v + D v - 2 \kappa (\tilde{u}_t+ u_t)\, v_t - 2\kappa \tilde{u}_{tt} \, v
\\
&= 2(\tilde{\kappa}-\kappa)(\tilde{u}\tilde{u}_{tt} + \tilde{u}_t^2)\mbox{ in }\Omega\times(0,T)\\
v(0)=0, \quad v_t(0)&=0 \mbox{ in }\Omega
\end{aligned}
\end{equation}
as well as the first order Taylor remainder $w=G(\tilde{\kappa})-G(\kappa)-G'(\kappa)(\tilde{\kappa}-\kappa)$ which satisfies
\begin{equation}\label{eqn:w}
\begin{aligned}
(1-2\kappa u)w_{tt}&+c^2\mathcal{A} w + D w - 4 \kappa u_t\, w_t - 2\kappa u_{tt} \, w\\
&= 2\underline{d\kappa}(v \tilde{u}_{tt} + u v_{tt} + (\tilde{u}_t+u_t) v_t ) + 2\kappa (v v_{tt} + v_t^2)
\mbox{ in }\Omega\times(0,T)\\
w(0)=0, \quad w_t(0)&=0 \mbox{ in }\Omega \,,
\end{aligned}
\end{equation}

Here we can allow for spatially varying sound speed $c_0(x)$ for which we only require 
\begin{equation}\label{eqn:c0}
c_0\in L^\infty(\Omega) \mbox{ and }
c_0(x)\geq c>0
\end{equation}
unless otherwise stated, by setting 
\begin{equation}\label{eqn:opA}
\mathcal{A}={{\small -\frac{c_0^2}{c^2}}}\,\triangle
\end{equation} 
where $-\triangle$ is the Laplace operator equipped with homogeneous Dirichlet boundary conditions.
We denote by $(\phi_j,\lambda_j)$ an eigensystem of the operator $\mathcal{A}$ with domain $\dot{H}^2(\Omega):=\mathcal{D}(\mathcal{A})$ which is selfadjoint and positive definite with respect to the weighted $L^2$ space $\dot{L}^2(\Omega):=L^2_{c^2/c_0^2}(\Omega)$. Note that we assume $\mathcal{A}^{-1}:\dot{L}^2(\Omega)\to\dot{L}^2(\Omega)$ to be compact (which is e.g., satisfied if $\Omega$ is bounded; for some comments on more general domain and boundary settings we point to \cite{KaltenbacherRundell:2021b}) , so that the eigensystem exists and is complete with $\lambda_j\to\infty$ as $j\to\infty$. Moreover, this defines a scale of Hilbert spaces $\dot{H}^s(\Omega):=\mathcal{D}(\mathcal{A}^{s/2})$, $s\in\mathbb{R}$, whose norm can be defined via the eigensystem as $\|v\|_{\dot{H}^s(\Omega)}=\left(\sum_{j=1}^\infty \lambda_j^s |\langle v,\phi_j\rangle|^2\right)^{1/2}$.
We will denote by $\langle\cdot,\cdot\rangle$ the $\dot{L}^2$ inner product (that is, the weighted one) on $\Omega$ whereas the use of the ordinary $L^2$ inner product will be indicated by a subscript $\langle\cdot,\cdot\rangle_{L^2}$
Moreover, we use the abbreviations   
$\|u\|_{L^p_t(L^q)}=\|u\|_{L^p(0,t;L^q(\Omega))}$, 
$\|u\|_{L^p(L^q)}=\|u\|_{L^p(0,T;L^q(\Omega))}$
for space-time norms.
\def\LoneLtwot{L^1_t(L^2)} 
\def\LtwoLtwot{L^2_t(L^2)} 
\def\LinfLtwot{L^\infty_t(L^2)}
\def\LtwoLtwodt{L^2_t(\dot{L}^2)}
\def\LinfLtwodt{L^\infty_t(\dot{L}^2)}
\def\LtwoLinft{L^2_t(L^\infty)}
\def\LinfLinft{L^\infty_t(L^\infty)}
\def\LtwoLfourt{L^2(L^4)}
\def\LinfLfourt{L^\infty(L^4)}
\def\LtwoLtwo{L^2(L^2)}
\def\LtwoLtwod{L^2(\dot{L}^2)}
\def\LinfLtwo{L^\infty(L^2)}
\def\LtwoLinf{L^2(L^\infty)}
\def\LinfLinf{L^\infty(L^\infty)}
\def\LtwoLfour{L^2(L^4)}
\def\LinfLfour{L^\infty(L^4)}

Throughout this paper, we denote by $\partial_t^\alpha$ the (partial) Caputo-Djrbashian fractional time derivative of order $\alpha\in(n-1,n)$ with $n\in\mathbb{N}$ by
$\partial_t^\alpha u  = I_t^{n-\alpha}[\partial^n_t u]\;$,
where $\partial^n_t$ denotes the $n$-th integer order partial time derivative and for
$\gamma\in(0,1)$, and $I_t^\gamma$ is the Abel integral operator defined by 
\[
I_t^\gamma[v](t) = \frac{1}{\Gamma(\gamma)}\int_{0}^{t}\frac{v(s)}{(t-s)^{1-\gamma}}ds\,.
\]
For details on fractional differentiation and subdiffusion equations,
we refer to, e.g.,
\cite{Dzjbashian:1966,Djrbashian:1993,
MainardiGorenflo:2000,SakamotoYamamoto:2011a,SamkoKilbasMarichev:1993}.
See also the tutorial paper on inverse problems for anomalous diffusion processes
\cite{JinRundell:2015}.
Whenever we use the Riemann-Liouville fractional derivative $\partial^n_t I_t^{n-\alpha}$, this will be denoted by ${}^{RL}\partial_t^\alpha$.
These two versions of the fractional derivative coincide when applied to a function whose Cauchy data up to order $n-1$ at $t=0$ vanish.

The crucial tool we need in obtaining the required estimates is the following consequence of Alikhanov's Lemma \cite[Lemma 1]{Alikhanov:11}
\begin{equation}\label{eqn:Alikhanov_1}
\partial_t^{\gamma} [w](t) w(t)\geq \tfrac12(\partial_t^{\gamma} w^2)(t)
\end{equation}
for any absolutely continuous function $w$. 
We apply it to $w=\partial_t^\alpha v$ with $\gamma=1-\alpha$, using the identities 
\[
\begin{aligned}
&\qquad \partial_t^{1-\alpha}w
= \partial_t^{1-\alpha}\partial_t^{\alpha} v
= \partial_t^{1-\alpha} I_t^{1-\alpha} v_t
= v_t
\,, \\
&\qquad\int_0^t (\partial_t^{1-\alpha} w^2)(s)\, ds = 
\int_0^t {}^{RL}\partial_t^{1-\alpha} [w^2](s)\, ds = 
\int_0^t \partial_t I_t^\alpha [w^2](s)\, ds 
= I_t^\alpha [w^2](t) 
\end{aligned}
\]
that hold for $v_t\in L^\infty(0,T)$ (Theorem 2.3 in fde-lect.pdf) 
and for $w^2\in W^{1,1}(0,T)$ with $w(0)=0$ (Theorem 2.2 in fde-lect.pdf).
Note that for $w= \partial_t^\alpha v$ we automatically have $w(0)=0$ and $I_t^\alpha [w^2](t)=0$.
After integration with respect to time this implies the following result.

\begin{lemma}\label{lem:Alikhanov2}
For $v\in W^{1,\infty}(0,T)$ with $(\partial_t^\alpha v)^2\in W^{1,1}(0,T)$, and $t\in(0,T)$, the following estimate holds.
\begin{equation}\label{eqn:Alikhanov_2}
\int_0^t\partial_t^{\alpha} [v](s) v_t(s)\, ds 
\geq I_t^\alpha\left[(\partial_t^{\alpha} v)^2\right]
\geq \frac{1}{2\Gamma(\alpha)t^{1-\alpha}}\|\partial_t^{\alpha} v\|_{L^2(0,t)}^2\,.
\end{equation}
\end{lemma}

\subsection{Caputo-Wismer-Kelvin damping}
We start with the Caputo-Wismer-Kelvin model 
\[
D = b \mathcal{A}^\beta \partial_t^\alpha  \quad\mbox{ with }\beta\in[0,1], \ b\geq0
\]
and first of all consider the initial boundary value problem for the general linear {\sc pde}
\begin{equation}\label{eqn:linPDE_CWK}
(1-\sigma)u_{tt}+c^2\mathcal{A} u + b \mathcal{A}^\beta \partial_t^\alpha u + \mu u_t + \rho u
= h
\end{equation}
\begin{equation}\label{eqn:init_CWK}
u(0)= u_0\,, \quad u_t(0)=u_1 
\end{equation}
with constants $b$, $c>0$ and given space- and time dependent functions $\sigma$, $\mu$, $\rho$, $h$ where $\sigma$ satisfies the non-degeneracy condition
\begin{equation}\label{eqn:nondegeneracy}
\sigma(x,t)\leq \overline{\sigma}<1 \mbox{ for all }x\in\Omega\, \quad t\in(0,T)\,.
\end{equation}

In order to prove existence and uniqueness of solutions to \eqref{eqn:linPDE_CWK}, \eqref{eqn:init_CWK},
we apply the usual Faedo-Galerkin approach of discretisation in space with eigenfunctions of $\mathcal{A}$, 
$u(x,t)\approx u^n(x,t)=\sum_{i=1}^n u^n_i(t)\phi_i(x)$ and testing with $\phi_j$, that is,
\begin{equation}\label{eqn:Galerkin_CWK0}
\langle (1-\sigma)u^n_{tt}+c^2\mathcal{A} u^n + b \mathcal{A}^\beta \partial_t^\alpha u^n + \mu u^n_t + \rho u^n
- h,v\rangle = 0 \, \quad v\in\mbox{span}(\phi_1,\ldots, \phi_n)\,.
\end{equation}
This leads to the {\sc ode} system
\begin{equation}\label{eqn:Galerkin_CWK}
(I-S^n(t)) {\underline{u}^n}''(t) 
+ b \, (\Lambda^n)^\beta (\partial_t^\alpha \underline{u}^n)(t)
+ M^n(t) {\underline{u}^n}'(t)
+ \Bigl(c^2 \, \Lambda^n +R^n(t)\Bigr) \underline{u}^n(t)
= \underline{h}^n
\end{equation}
with matrices and vectors defined by 
\begin{equation}\label{eqn:matricesGalerkin_CWK}
\begin{aligned}
&\underline{u}^n (t)= (u^n_i(t))_{i=1,\ldots n}\,, \quad 
\underline{h}^n (t) = (\langle h(t),\phi_i\rangle)_{i=1,\ldots n}\,, \quad
\Lambda^n = \mbox{diag}(\lambda_1,\ldots,\lambda_n)\,,
\\&
S^n(t)=(\langle\sigma(t)\phi_i,\phi_j\rangle)_{i,j=1,\ldots n}\,, \quad 
M^n(t)=(\langle\mu(t)\phi_i,\phi_j\rangle)_{i,j=1,\ldots n}\,, \\\
&R^n(t)=(\langle\rho(t)\phi_i,\phi_j\rangle)_{i,j=1,\ldots n}\,.
\end{aligned}
\end{equation}
Existence of a unique solution $\underline{u}^n\in C^2(0,T;\mathbb{R}^n)$ to \eqref{eqn:Galerkin_CWK} follows from standard {\sc ode} theory (Picard-Lindel\"of Theorem and Gronwall's Inequality), as long as $\sigma$, $\mu$ and $\rho$ are in $C(0,T;\dot{H}^s(\Omega))$ for some $s\in\mathbb{R}$ (noting that the eigenfunctions $\phi_j$ are contained in $\dot{H}^k(\Omega)$ for any $k\in\mathbb{N}$ and therefore the vector and matrix functions $\underline{h}^n$, $S^n$, $M^n$, $R^n$ in \eqref{eqn:matricesGalerkin_CWK} are well-defined and contained in $C(0,T;\mathbb{R}^n)$ and $C(0,T;\mathbb{R}^{n\times n})$, respectively. 
Moreover due to \eqref{eqn:nondegeneracy}, the symmetric matrix $S^n(t)$ is positive definite with smallest eigenvalue bounded away from zero by $1-\overline{\sigma}$ cf \eqref{eqn:nondegeneracy}.

We multiply \eqref{eqn:Galerkin_CWK} with $(\Lambda^n)^2 {\underline{u}^n}'(t)$ and integrate with respect to time, using the identity 
\[
\begin{aligned}
&(S^n(t) {\underline{u}^n}''(t))^T (\Lambda^n)^2 {\underline{u}^n}'(t)
= \sum_{i=1}^n \sum_{j=1}^n \langle\sigma(t)\phi_i,\phi_j\rangle \lambda_j^2 {u^n_i}''(t) {u^n_j}'(t)\\
&= \sum_{i=1}^n \sum_{j=1}^n \langle\sigma(t)\phi_i,\phi_j\rangle \lambda_i {u^n_i}''(t) \lambda_j {u^n_j}'(t)
+ \sum_{i=1}^n \sum_{j=1}^n \langle\sigma(t)\phi_i,\phi_j\rangle (\lambda_j-\lambda_i) {u^n_i}''(t) \lambda_j {u^n_j}'(t)
\\
&= \frac12 \frac{d}{dt} \sum_{i=1}^n \sum_{j=1}^n \langle\sigma(t)\phi_i,\phi_j\rangle \lambda_i {u^n_i}'(t) \lambda_j {u^n_j}'(t) 
- \frac12 \sum_{i=1}^n \sum_{j=1}^n \langle\sigma_t(t)\phi_i,\phi_j\rangle \lambda_i {u^n_i}'(t) \lambda_j {u^n_j}'(t)\\
&\qquad + \sum_{i=1}^n \sum_{j=1}^n \Bigl(\langle\sigma(t)\phi_i,\mathcal{A}\phi_j\rangle-\langle\sigma(t)\phi_j,\mathcal{A}\phi_i\rangle\Bigr)  {u^n_i}''(t) \lambda_j {u^n_j}'(t)\\
\end{aligned}
\]
where
\[
\begin{aligned}
&\langle\sigma(t)\phi_i,\mathcal{A}\phi_j\rangle-\langle\sigma(t)\phi_j,\mathcal{A}\phi_i\rangle
= \langle\mathcal{A}[\sigma(t)\phi_i]-\sigma(t)\mathcal{A}\phi_i,\phi_j\rangle\\
&= \langle-\triangle[\sigma(t)\phi_i]+\sigma(t)\triangle\phi_i,\phi_j\rangle_{L^2}
= -\langle(c_0^2/c^2)(\triangle \sigma(t)\, \phi_i+\nabla\sigma(t)\cdot\nabla\phi_i),\phi_j\rangle
\end{aligned}
\]
provided $\sigma(t)\phi_i\in \dot{H}^s(\Omega))$ for some $s\in\mathbb{R}$.
(Note that the latter identity also holds true in case of spatially varying $c_0$ since we use the $c^2/c_0^2$ weighted $L^2$ inner product then).
Thus we have 
\[
\begin{aligned}
&\int_0^t((I-S^n(t)) {\underline{u}^n}''(s))^T (\Lambda^n)^2 {\underline{u}^n}'(s)\, ds\\
&= \int_0^t\Bigl(\frac12 \frac{d}{dt} \|\sqrt{1-\sigma}\mathcal{A}u^n_t\|_{\dot{L}^2(\Omega)}^2(t)
+ \frac12 \langle\sigma_t(t)\mathcal{A}u^n_t(t),\mathcal{A}u^n_t(t)\rangle\\
&\qquad+ \langle(c_0^2/c^2)(\triangle \sigma(t)\, u^n_{tt}(t)+\nabla\sigma(t)\cdot\nabla u^n_{tt}(t)),\mathcal{A}u^n_t(t)\rangle
\Bigr)\\
&= \frac12 \|\sqrt{1-\sigma(t)}\mathcal{A}u^n_t(t)\|_{\dot{L}^2(\Omega)}^2- \frac12 \|\sqrt{1-\sigma(0)}\mathcal{A}u^n_t(0)\|_{\dot{L}^2(\Omega)}^2\\
&\quad+\int_0^t\langle\frac12\sigma_t(s)\mathcal{A}u^n_t(s)
+(c_0^2/c^2)(\triangle\sigma(s)\, u^n_{tt}(s)+\nabla\sigma(s)\cdot\nabla u^n_{tt}(s)),\mathcal{A}u^n_t(s)\rangle\, ds
\Bigr)\,.
\end{aligned}
\]
Similarly, we have 
\[
\begin{aligned}
&(M^n(t) {\underline{u}^n}'(t))^T (\Lambda^n)^2 {\underline{u}^n}'(t)\\
&= \langle\mu(t)\mathcal{A}u^n_t(t),\mathcal{A}u^n_t(t)\rangle
- \langle(c_0^2/c^2)(\triangle \mu(t)\, u^n_t(t)+\nabla\mu(t)\cdot\nabla u^n_t(t)),\mathcal{A}u^n_t(t)\rangle
\end{aligned}
\]
and 
\[
\begin{aligned}
&(R^n(t) {\underline{u}^n}(t))^T (\Lambda^n)^2 {\underline{u}^n}'(t)\\
&= \langle\rho(t)\mathcal{A}u^n(t),\mathcal{A}u^n_t(t)\rangle
- \langle(c_0^2/c^2)(\triangle \rho(t)\, u^n(t)+\nabla\rho(t)\cdot\nabla u^n(t)),\mathcal{A}u^n_t(t)\rangle\,.
\end{aligned}
\]
Finally,
\[
\begin{aligned}
&\int_0^t((\Lambda^n)^\beta (\partial_t^\alpha \underline{u}^n)(t))^T(\Lambda^n)^2 {\underline{u}^n}'(t)
=\int_0^t\sum_{j=1}^n \lambda_j^{2+\beta} \int_0^t(\partial_t^\alpha u^n_j)(s) {u^n_j}'(s)\, ds\\
&\geq \frac{1}{2\Gamma(\alpha)t^{1-\alpha}} \sum_{j=1}^n \lambda_j^{2+\beta} \int_0^t\Bigl(\partial_t^{\alpha} u^n_j(s)\Bigr)^2\, ds
=\frac{1}{2\Gamma(\alpha)t^{1-\alpha}} \|\mathcal{A}^{1+\beta/2}\partial_t^{\alpha}u^n\|_{\LtwoLtwodt}^2
\,.
\end{aligned}
\]
and 
\[
\begin{aligned}
&\int_0^t(\Lambda^n \underline{u}^n(t))^T(\Lambda^n)^2 {\underline{u}^n}'(t)
=\sum_{j=1}^n \lambda_j^3 \int_0^t u^n_j(s) {u^n_j}'(s)\, ds\\
&=\frac12 \sum_{j=1}^n \lambda_j^3 \int_0^t \frac{d}{dt}\Bigl(u^n_j\Bigr)^2(s)\, ds
=\frac12 \|\mathcal{A}^{3/2} u^n(t)\|_{\dot{L}^2(\Omega)}^2 - \frac12 \|\mathcal{A}^{3/2} u^n(0)\|_{\dot{L}^2(\Omega)}^2 
\,.
\end{aligned}
\]
This together with Young's inequality yields the energy estimate
\begin{equation}\label{eqn:enest1}
\begin{aligned}
&\frac12 \|\sqrt{1-\sigma(t)}\mathcal{A}u^n_t(t)\|_{\dot{L}^2(\Omega)}^2
+\frac{b}{2\Gamma(\alpha)t^{1-\alpha}} \|\mathcal{A}^{1+\beta/2}\partial_t^{\alpha}u^n\|_{\LtwoLtwodt}^2
+\frac{c^2}{2} \|\mathcal{A}^{3/2} u^n(t)\|_{\dot{L}^2(\Omega)}^2 
\\
&\leq\frac12 \|\sqrt{1-\sigma(0)}\mathcal{A}u^n_t(0)\|_{\dot{L}^2(\Omega)}^2
+\frac{c^2}{2} \|\mathcal{A}^{3/2} u^n(0)\|_{\dot{L}^2(\Omega)}^2 
+ \frac{1}{2\epsilon} \|\mathcal{A}u^n_t\|_{\LtwoLtwodt}^2\\
&\ + \frac{\epsilon}{2} \|-\frac12\sigma_t\mathcal{A}u^n_t - \mu\mathcal{A}u^n_t -\rho\mathcal{A}u^n
+\mathcal{A}h\\
&\quad\ 
+(c_0^2/c^2)(
-\triangle\sigma\, u^n_{tt}-\nabla\sigma\cdot\nabla u^n_{tt} +\triangle \mu\, u^n_t+\nabla\mu\cdot\nabla u^n_t
+\triangle \rho\, u^n+\nabla\rho\cdot\nabla u^n)\|_{\LtwoLtwodt}^2\\
&\leq\frac12 \|\sqrt{1-\sigma(0)}\mathcal{A}u^n_t(0)\|_{\dot{L}^2(\Omega)}^2
+\frac{c^2}{2} \|\mathcal{A}^{3/2} u^n(0)\|_{\dot{L}^2(\Omega)}^2 
+ \frac{1}{2\epsilon} \|\mathcal{A}u^n_t\|_{\LtwoLtwodt}^2\\
&\quad + 
\frac{\epsilon}{2} \Bigl(\frac12\|\sigma_t\|_{\LtwoLinft}\|\mathcal{A}u^n_t\|_{\LinfLtwot}
+\|\mu\|_{\LtwoLinft} \|\mathcal{A}u^n_t\|_{\LinfLtwot}\\
&\qquad\qquad
+\|\rho\|_{\LtwoLinft} \|\mathcal{A}u^n\|_{\LinfLtwot}
+ \|\mathcal{A}h\|_{\LtwoLtwot}\\
&\qquad\quad+\frac{\|c_0\|_{L^\infty(\Omega)}}{c} \bigl(
\|\triangle\sigma\|_{\LinfLfourt} \|u^n_{tt}\|_{\LtwoLfourt} 
+\|\nabla\sigma\|_{\LinfLinft} \|\nabla u^n_{tt}\|_{\LtwoLtwot}
\\
&\qquad \qquad +\|\triangle \mu\|_{\LtwoLtwot}\|u^n_t\|_{\LinfLinft}
+\|\nabla\mu\|_{\LtwoLfourt} \|\nabla u^n_t\|_{\LinfLfourt}\\
&\qquad \qquad 
+\|\triangle \rho\|_{\LtwoLtwot} \|u^n\|_{\LinfLinft} 
+\|\nabla\rho\|_{\LtwoLfourt} \|\nabla u^n\|_{\LinfLfourt}
\bigl)\Bigr)^2
\end{aligned}
\end{equation}
Here we can make use of the fact that
$\|u^n\|_{L^\infty(0,t;Z)}\leq \|u^n(0)\|_{Z} + \sqrt{T} \|u^n_t\|_{L^2(0,t;Z)}$
and the embedding estimates
\begin{equation}\label{eqn:embeddings}
\begin{aligned}
&\|v\|_{L^4(\Omega)}\leq C_{H^1,L^4}\|\nabla v\|\,, \quad v\in H_0^1(\Omega)\,, \\
&\|v\|_{L^\infty(\Omega)}\leq C_{H^2,L^\infty}\|\mathcal{A} v\|_{\dot{L}^2(\Omega)}\,, \quad v\in H_0^1(\Omega)\cap H^2(\Omega)\,,
\end{aligned}
\end{equation}
in order to further estimate
\[
\begin{aligned}
\|u^n_{tt}\|_{\LtwoLfourt}&\leq C_{H^1,L^4} \|\nabla u^n_{tt}\|_{\LtwoLtwot}\\
\|u^n_t\|_{\LinfLinft}&\leq C_{H^2,L^\infty} \|\mathcal{A}u^n_t\|_{\LinfLtwot}\\
\|u^n\|_{\LinfLinft}&\leq C_{H^2,L^\infty} \|\mathcal{A}u^n\|_{\LinfLtwot}\\
\|\nabla u^n\|_{\LinfLinft}&\leq C_{H^2,L^\infty} \|\mathcal{A}^{3/2}u^n\|_{\LinfLtwot}
\end{aligned}
\]
Now we proceed with estimating $\|\nabla u^n_{tt}\|_{L^2}^2$ 
by multiplying the {\sc ode} \eqref{eqn:Galerkin_CWK} with $\Lambda^n {\underline{u}^n}''$, 
that is testing \eqref{eqn:Galerkin_CWK0} with $v=\mathcal{A} u^n_{tt}(t)$,
and using integration by parts (note that all terms in $(1-\sigma)u^n_{tt}+c^2\mathcal{A} u^n + b \mathcal{A}^\beta \partial_t^\alpha u^n + \mu u^n_t + \rho u^n-h$ vanish on $\partial\Omega$) as well as Young's inequality
\[
\begin{aligned}
&0=\langle (1-\sigma)u^n_{tt}+c^2\mathcal{A} u^n + b \mathcal{A}^\beta \partial_t^\alpha u^n + \mu u^n_t + \rho u^n
- h, \mathcal{A} u^n_{tt}\rangle\\
&=\|\sqrt{1\!-\!\sigma}\nabla u^n_{tt}\|_{L^2(\Omega)}^2
+ \langle - \nabla\sigma\, u^n_{tt} +\nabla\Bigl(c^2\mathcal{A} u^n + b \mathcal{A}^\beta \partial_t^\alpha u^n 
+ \mu u^n_t + \rho u^n \!-\! h\Bigr), 
\nabla u^n_{tt}\rangle_{L^2(\Omega)}\\
&\geq \|\sqrt{1-\sigma}\nabla u^n_{tt}\|_{L^2(\Omega)}^2-\frac12 (1-\overline{\sigma})\|\nabla u^n_{tt}\|_{L^2(\Omega)}^2\\
&\quad-\frac{1}{2(1-\overline{\sigma})}
\|-\nabla \sigma\, u^n_{tt} +c^2\nabla\mathcal{A} u^n + b \nabla\mathcal{A}^\beta \partial_t^\alpha u^n \\
&\hspace*{3cm}+ \nabla\mu \, u^n_t + \mu \nabla u^n_t + \nabla\rho \, u^n + \rho \nabla u^n -\nabla h\|_{L^2(\Omega)}^2
\end{aligned}
\]
which yields 
\[
\begin{aligned}
&\|\nabla u^n_{tt}\|_{\LtwoLtwot}\leq 
\frac{1}{1-\overline{\sigma}}\\
&\cdot\Bigl(
\|\nabla \sigma\|_{\LinfLfourt} \|u^n_{tt}\|_{\LtwoLfourt}^2
+c^2 \|\mathcal{A}^{3/2} u^n\|_{\LtwoLtwot}
+b \|\mathcal{A}^{1/2+\beta} \partial_t^\alpha u^n\|_{\LtwoLtwot}\\
&\qquad+\|\nabla\mu\|_{\LtwoLfourt} \|u^n_t\|_{\LinfLfourt}
+\|\mu\|_{\LtwoLinft} \|\nabla u^n_t\|_{\LinfLtwot}\\
&\qquad
+\|\nabla\rho\|_{\LtwoLtwot} \|u^n\|_{\LinfLinft}
+\|\rho\|_{\LtwoLinft} \|\nabla u^n\|_{\LinfLtwot}\Bigr)
\end{aligned}
\]
where we can again employ the embedding estimates \eqref{eqn:embeddings} and assume 
\begin{equation}\label{eqn:nablasigmasmall}
\|\nabla \sigma\|_{\LinfLfourt} < \frac{1-\overline{\sigma}}{C_{H^1\to L^4}}
\end{equation}
in order to extract an estimate of the form
\begin{equation}\label{eqn:enest2}
\begin{aligned}
&\|\nabla u^n_{tt}\|_{\LtwoLtwot}\\
&\leq C \Bigl( 
\|\mathcal{A}^{3/2} u^n\|_{\LtwoLtwot} 
+ \|\mathcal{A}^{1/2+\beta} \partial_t^\alpha u^n\|_{\LtwoLtwot}
+ \|\nabla u^n_t\|_{\LinfLtwot}
+ \|\nabla u^n\|_{\LinfLtwot}
\Bigr)
\end{aligned}
\end{equation}

Adding a multiple (factor $\epsilon (\|\triangle\sigma\|_{\LinfLfourt}^2 (C_{H^1,L^4})^2 +\|\nabla\sigma\|_{\LinfLinft}^2)$) of the square of \eqref{eqn:enest2} to \eqref{eqn:enest1}, making $\epsilon$ small enough (so that $C^2\epsilon < \frac{b}{2\Gamma(\alpha)t^{1-\alpha}})$ and all terms containing $L^\infty(0,T)$ norms of $u^n$ on the right hand side of \eqref{eqn:enest1} can be dominated by left hand side terms) using the fact that $1/2+\beta\leq 1+\beta/2$ for $\beta\in[0,1]$ and  Gronwall's inequality we end up with an estimate of the form
\begin{equation}\label{eqn:est_linCWKn}
\begin{aligned}
&\|\nabla u^n_{tt}\|_{\LtwoLtwo}^2 
+ \|\mathcal{A}u^n_t\|_{\LinfLtwodt}^2
+\|\mathcal{A}^{1+\beta/2}\partial_t^{\alpha}u^n\|_{\LtwoLtwodt}^2
+\|\mathcal{A}^{3/2} u^n\|_{\LinfLtwodt}^2\\
&\leq C(T) \Bigl(\|\mathcal{A}u^n_t(0)\|_{\dot{L}^2(\Omega)}^2
+\|\mathcal{A}^{3/2} u^n(0)\|_{\dot{L}^2(\Omega)}^2
+\|\mathcal{A}h\|_{\LtwoLtwodt}^2\Bigr)
\end{aligned}
\end{equation}
which via weak limits shows existence of a solution to the homogeneous initial boundary value problem for \eqref{eqn:linPDE_CWK} and transfers to $u$ as
\begin{equation}\label{eqn:est_linCWK}
\begin{aligned}
&\|u\|_U^2:=\|\nabla u_{tt}\|_{\LtwoLtwo}^2 
+ \|\mathcal{A}u_t\|_{\LinfLtwodt}^2
+\|\mathcal{A}^{1+\beta/2}\partial_t^{\alpha}u\|_{\LtwoLtwodt}^2
+\|\mathcal{A}^{3/2} u\|_{\LinfLtwodt}^2\\
&\qquad\leq C(T) \Bigl(\|\mathcal{A}u_t(0)\|_{\dot{L}^2(\Omega)}^2
+\|\mathcal{A}^{3/2} u(0)\|_{\dot{L}^2(\Omega)}^2
+\|\mathcal{A}h\|_{\LtwoLtwodt}^2\Bigr)
\end{aligned}
\end{equation}

The required regularity on $\sigma$, $\mu$, $\rho$, $h$, $u_0$, $u_1$ is, besides \eqref{eqn:nondegeneracy} and \eqref{eqn:nablasigmasmall}
\begin{equation}\label{eqn:regularity_sigmamurhoh}
\begin{aligned}
&\sigma\in H^1(0,T;L^\infty(\Omega))\cap L^\infty(0,T;W^{2,4}(\Omega)\cap W^{1,\infty}(\Omega))\\
&\mu\in L^2(0,T;H^2(\Omega))\,, \quad \rho\in L^2(0,T;H^2(\Omega))
\\
&h\in L^2(0,T;\dot{H}^2(\Omega))\,, \quad
u_0\in \dot{H}^3(\Omega)\,, \quad
u_1\in \dot{H}^2(\Omega)
\end{aligned}
\end{equation}
\begin{proposition} \label{prop:wellposed_lin_CWK}
Under conditions \eqref{eqn:nondegeneracy}, \eqref{eqn:nablasigmasmall}, \eqref{eqn:regularity_sigmamurhoh}, there exists a unique solution 
\begin{equation}\label{defU}
u\in U:=H^2(0,T;\dot{L}^2(\Omega))\cap W^{1,\infty}(0,T;\dot{H}^2(\Omega))\cap L^\infty(0,T;\dot{H}^3(\Omega))
\end{equation}
to the initial boundary value problem 
\eqref{eqn:linPDE_CWK}, \eqref{eqn:init_CWK}. This solution satisfies the estimate \eqref{eqn:est_linCWK}.
\end{proposition}

The energy estimate leading to this result has been obtained by basically ``multiplying \eqref{eqn:linPDE_CWK} with $\mathcal{A}^2u_t$'', that is, taking the $\dot{L}^2(\Omega)$ inner product of the {\sc pde} with $\mathcal{A}^2u_t$ and using selfadjointness of $\mathcal{A}$ in $\dot{L}^2(\Omega)$.

Later on, we will also need less regular solutions along with estimates on them. Since the proofs are actually somewhat simpler then, we skip the details on Galerkin approximation and only provide the energy estimates.

Multiplying \eqref{eqn:linPDE_CWK} with $\mathcal{A}u_t$ we obtain
\begin{equation}\label{eqn:est_linCWK_low0}
\begin{aligned}
&\frac12 \|\sqrt{1-\sigma(t)}\nabla u_t(t)\|_{L^2(\Omega)}^2
+\frac{c^2}{2} \|\mathcal{A} u(t)\|_{\dot{L}^2(\Omega)}^2
+\frac{1}{2\Gamma(\alpha)t^{1-\alpha}}\|\mathcal{A}^{(1+\beta)/2} \partial_t^\alpha u\|_{\LtwoLtwodt}^2\\
&\leq \frac12 \|\sqrt{1-\sigma(0)}\nabla u_t(0)\|_{L^2(\Omega)}^2
+\frac{c^2}{2} \|\mathcal{A} u(0)\|_{\dot{L}^2(\Omega)}^2
\\
&\quad-\int_0^t\langle\frac12\sigma_t(s)\nabla u_t(s)+\nabla\sigma(s)\,u_{tt}(s)
+\nabla\mu\, u_t+\mu\nabla u_t+\nabla\rho\, u+\rho\nabla u\!-\!\nabla h,\nabla u_t(s)\rangle_{L^2(\Omega)}\, ds\\
&\leq \frac12 \|\sqrt{1-\sigma(0)}\nabla u_t(0)\|_{L^2(\Omega)}^2
+\frac{c^2}{2} \|\mathcal{A} u(0)\|_{\dot{L}^2(\Omega)}^2
+\frac{1}{2\epsilon} \|\nabla u_t\|_{\LtwoLtwot}^2\\
&\qquad +\frac{\epsilon}{2}\Bigl(
\frac12\|\sigma_t\|_{\LtwoLinft}\|\nabla u_t\|_{\LinfLtwot}
+\|\nabla\sigma\|_{\LinfLinft}\|u_{tt}\|_{\LtwoLtwot}
\\&\qquad\qquad
+\|\nabla\mu\|_{\LtwoLfourt}\|u_t\|_{\LinfLfourt}
+\|\mu\|_{\LtwoLinft}\|\nabla u_t\|_{\LinfLtwot}
\\&\qquad\qquad
+\|\nabla\rho\|_{\LtwoLtwot}\|u\|_{\LinfLinft}
+\|\rho\|_{\LtwoLfourt}\|\nabla u\|_{\LinfLfourt}
+\|\nabla h\|_{\LtwoLtwot}\Bigr)^2\,.
\end{aligned}
\end{equation}
The {\sc pde} provides us with 
\[
\begin{aligned}
&\|u_{tt}\|_{\LtwoLtwot}=\|-\frac{1}{1-\sigma}\Bigl( c^2\mathcal{A} u + b \mathcal{A}^\beta \partial_t^\alpha u + \mu u_t + \rho u - h\Bigr)\|_{\LtwoLtwo}\\ 
&\leq C (\|\mu\|_{\LtwoLinft}+\|\rho\|_{\LtwoLfourt})\times\\
&\qquad\sup_{t\in(0,T)}\Bigl(\|\sqrt{1\!-\!\sigma(t)}\nabla u_t(t)\|_{L^2(\Omega)}^2
+\frac{c^2}{2} \|\mathcal{A} u(t)\|_{\dot{L}^2(\Omega)}^2
+\frac{1}{2\Gamma(\alpha)t^{1-\alpha}}\|\mathcal{A}^\beta \partial_t^\alpha u\|_{\LtwoLtwodt}^2
\Bigr)
\end{aligned}
\]
which, due to the condition $\beta\leq(1+\beta)/2$,
 allows us to dominate the $\epsilon u_{tt}$ term on the right hand side of \eqref{eqn:est_linCWK_low0}.
Thus, using Gronwall's inequality, we get an estimate of the form
\begin{equation}\label{eqn:est_linCWK_low}
\begin{aligned}
&\|u\|_{U_{lo}}^2:= \|u_{tt}\|_{\LtwoLtwodt}^2
+\|\nabla u_t\|_{\LinfLtwodt}^2
+\|\mathcal{A}^{(1+\beta)/2}\partial_t^{\alpha}u\|_{\LtwoLtwodt}^2
+\|\mathcal{A} u\|_{\LinfLtwodt}^2\\
&\leq C(T) \Bigl(\|\nabla u_t(0)\|_{L^2(\Omega)}^2
+\|\mathcal{A} u(0)\|_{\dot{L}^2(\Omega)}^2
+\|\nabla h\|_{\LtwoLtwo}^2\Bigr)
\end{aligned}
\end{equation}
provided
\begin{equation}\label{eqn:regularity_sigmamurhoh_low}
\begin{aligned}
&\sigma\in H^1(0,T;L^\infty(\Omega))\cap L^\infty(0,T;W^{1,\infty}(\Omega))\,, \\
&\mu\in L^2(0,T;L^\infty(\Omega)\cap W^{1,4}(\Omega))\,, \quad \rho\in L^2(0,T;H^1(\Omega))
\\
&h\in L^2(0,T;H^1(\Omega))\,, \quad
u_0\in \dot{H}^2(\Omega)\,, \quad
u_1\in \dot{H}^1(\Omega)
\end{aligned}
\end{equation}
\begin{proposition} \label{prop:wellposed_lin_CWK_low}
Under the conditions \eqref{eqn:nondegeneracy},
\eqref{eqn:regularity_sigmamurhoh_low}, there exists a unique solution $u$,
 where
 $u\in H^2(0,T;\dot{L}^2(\Omega))\cap W^{1,\infty}(0,T;\dot{H}^1(\Omega))\cap L^\infty(0,T;\dot{H}^2(\Omega))$,
 to the initial boundary value problem 
\eqref{eqn:linPDE_CWK}, \eqref{eqn:init_CWK}, and this solution satisfies the
estimate \eqref{eqn:est_linCWK_low}.
\end{proposition}

An even lower regularity estimate can be obtained by multiplication of \eqref{eqn:linPDE_CWK} with $u_t$, which yields
\[
\begin{aligned}
&\frac12 \|\sqrt{1-\sigma(t)}u_t(t)\|_{\dot{L}^2(\Omega)}^2
+\frac{c^2}{2} \|\nabla u(t)\|_{L^2(\Omega)}^2
+\frac{1}{2\Gamma(\alpha)t^{1-\alpha}}\|\mathcal{A}^{\beta/2} \partial_t^\alpha u\|_{\LtwoLtwodt}^2\\
&\leq \frac12\|\sqrt{1-\sigma(0)}u_t(0)\|_{\dot{L}^2(\Omega)}^2+\frac{c^2}{2} \|\nabla u(0)\|_{L^2(\Omega)}^2 \\
&\qquad- \int_0^t\langle\frac12\sigma_t(s)u_t(s) +\mu\, u_t+\rho u- h,u_t(s)\rangle\, ds\\
&\leq \frac12\|\sqrt{1-\sigma(0)}u_t(0)\|_{\dot{L}^2(\Omega)}^2+\frac{c^2}{2} \|\nabla u(0)\|_{L^2(\Omega)}^2  
+ \frac{1}{2\epsilon} \|u_t\|_{\LtwoLtwodt}^2\\
&\qquad
+\frac{\epsilon}{2} \Bigl(
\frac12\|\sigma_t\|_{\LtwoLinft}\|u_t\|_{\LinfLtwodt}
+\|\mu\|_{\LtwoLinft}\|u_t\|_{\LinfLtwodt}\\
&\qquad\qquad+\|\rho\|_{\LtwoLfourt}\|u\|_{\LinfLfourt}
+\|h\|_{\LtwoLtwodt}\Bigr)^2\,, 
\end{aligned}
\]
(where we have used the fact that $\|\frac{c}{c_0}\|_{L^\infty(\Omega)}\leq1$),
hence an estimate of the form
\begin{equation}\label{eqn:est_linCWK_verylow}
\begin{aligned}
&\|u\|_{U_{vl}}^2:= \ \|u_t\|_{\LinfLtwodt}^2
+\|\mathcal{A}^{\beta/2}\partial_t^{\alpha}u\|_{\LtwoLtwo}^2
+\|\nabla u\|_{\LinfLtwot}^2\\
&\leq C(T) \Bigl(\|u_t(0)\|_{\dot{L}^2(\Omega)}^2
+\|\nabla u(0)\|_{L^2(\Omega)}^2
+\|h\|_{\LtwoLtwodt}^2\Bigr)
\end{aligned}
\end{equation}
Here it obviously suffices to assume 
\begin{equation}\label{eqn:regularity_sigmamurhoh_verylow}
\begin{aligned}
&\sigma\in H^1(0,T;L^\infty(\Omega))\,, \quad 
\mu\in \LtwoLinft\,, \quad \rho\in \LtwoLfourt
\\
&h\in \LtwoLtwodt\,, \quad
u_0\in \dot{H}^1(\Omega)\,, \quad
u_1\in \dot{L}^2(\Omega)
\end{aligned}
\end{equation}
\begin{proposition} \label{prop:wellposed_lin_CWK_verylow}
Under conditions \eqref{eqn:nondegeneracy} and \eqref{eqn:regularity_sigmamurhoh_verylow}, there exists a unique solution $u$ lying in the space
 $u\in H^2(0,T;\dot{H}^{-1}(\Omega))\cap W^{1,\infty}(0,T;\dot{L}^2(\Omega))\cap L^\infty(0,T;\dot{H}^1(\Omega))$
and satisfying the initial boundary value problem 
\eqref{eqn:linPDE_CWK}, \eqref{eqn:init_CWK}.
This solution satisfies the estimate \eqref{eqn:est_linCWK_verylow}.
\end{proposition}

\medskip
We proceed to proving well-posedness of the nonlinear problem \eqref{eqn:Westervelt_init_D} with 
$D = b \mathcal{A}^\beta \partial_t^\alpha$ by applying a fixed point argument to the operator $\mathcal{T}$ mapping $v$ to the solution of 
\begin{equation}\label{eqn:Westervelt_init_D_fixedpoint}
\begin{aligned}
(1-2\kappa v)u_{tt}+c^2\mathcal{A} u + b \mathcal{A}^\beta\partial_t^\alpha u -2\kappa v_tu_t&=r \mbox{ in }\Omega\times(0,T)\\
u(0)=0, \quad u_t(0)&=0 \mbox{ in }\Omega
\end{aligned}
\end{equation}
that is, of \eqref{eqn:linPDE_CWK}, \eqref{eqn:init_CWK} with $\sigma(x,t)=-2\kappa(x)v(x,t)$, $\mu(x,t)=-2\kappa(x)v_t(x,t)$, $\rho(x,t)=0$, $h(x,t)=r(x,t)$.
We first consider self-mapping of $\mathcal{T}$.
Even in case of constant $\kappa$, the regularity requirements on $\sigma$, $\mu$, $\rho$ force us into the high regularity scenario of Proposition~\oldref{prop:wellposed_lin_CWK}. 
For spatially variable $\kappa$, due to the estimates 
\[
\begin{aligned}
\|\triangle&\mu\|_{\LtwoLtwo}= 2\|\triangle(\kappa v_t)\|_{\LtwoLtwo} \\
\quad&\leq 2 \Bigl(\|\triangle\kappa\|_{L^2(\Omega)}\|v_t\|_{L^2(0,T;L^\infty(\Omega)}
+2\|\nabla\kappa\|_{L^4(\Omega)}\|\nabla v_t\|_{\LtwoLfour}
+\|\kappa\|_{L^\infty(\Omega)}\|\triangle v_t\|_{\LtwoLtwo}
\Bigr)\\
\|\sigma_t&\|_{\LtwoLinft}\leq 2C_{H^2,L^\infty} \|\triangle(\kappa v_t)\|_{\LtwoLtwo}\\
\end{aligned}
\]
\[
\begin{aligned}
\|\triangle\sigma&\|_{\LinfLfourt}= 2\|\triangle(\kappa v)\|_{\LinfLfourt} \\
&\quad\leq 2 \Bigl(\|\triangle\kappa\|_{L^4(\Omega)}\|v\|_{\LinfLinf}
+2\|\nabla\kappa\|_{L^4(\Omega)}\|\nabla v\|_{\LinfLinf}
+\|\kappa\|_{L^\infty(\Omega)}\|\triangle v\|_{\LinfLfour}
\Bigr)\\
\|\nabla\sigma&\|_{\LinfLinft}= 2\|\nabla(\kappa v)\|_{\LinfLinft} 
\leq 2 \Bigl(\|\nabla\kappa\|_{L^\infty(\Omega)}\|v\|_{\LinfLinf}
+\|\kappa\|_{L^\infty(\Omega)}\|\nabla v\|_{\LinfLinf}
\Bigr)
\end{aligned}
\]
the regularity 
\begin{equation}\label{eqn:regkappa}
\kappa \in W^{2,4}(\Omega)\cap W^{1,\infty}(\Omega)
\end{equation}
is sufficient for obtaining the regularity \eqref{eqn:regularity_sigmamurhoh} for any for any $v\in U$.
To achieve the nondegeneracy and smallness conditions \eqref{eqn:nondegeneracy}, \eqref{eqn:nablasigmasmall}, we use the estimates
\begin{equation}\label{eqn:smallness_v}
\begin{aligned}
&\|\sigma\|_{\LinfLinft}\leq 2 \|\kappa\|_{L^\infty(\Omega)}\|v\|_{\LinfLinft}\\
&\|\nabla\sigma\|_{\LinfLfourt}
\leq 2 \|\nabla\kappa\|_{L^4(\Omega)}\|v\|_{\LinfLinf}
+2\|\kappa\|_{L^\infty(\Omega)}\|\nabla v\|_{\LinfLfour}
\end{aligned}
\end{equation}
and additionally to \eqref{eqn:regkappa} require smallness of $v$.

Proposition~\oldref{prop:wellposed_lin_CWK} yields that $\mathcal{T}$ is a self-mapping on 
\[
B_R = \{ v\in U \ : \ \|v\|_{U}\leq R\}
\]
provided the initial and right hand side data are sufficiently small so that 
\begin{equation}\label{eqn:smallness_init}
C(T) \Bigl(\|\mathcal{A}u_1\|_{\dot{L}^2(\Omega)}^2
+\|\mathcal{A}^{3/2} u_0\|_{\dot{L}^2(\Omega)}^2
+\|\mathcal{A}r\|_{\LtwoLtwodt}^2\Bigr)
\leq R^2
\end{equation}
with $C(T)$ as in Proposition~\oldref{prop:wellposed_lin_CWK}.
In view of \eqref{eqn:nondegeneracy}, \eqref{eqn:nablasigmasmall}, \eqref{eqn:smallness_v}, we choose $R$ such that 
\begin{equation}\label{eqn:smallness_R}
\Bigl(C_{H^1\to L^4}\Bigl(\|\nabla\kappa\|_{L^4(\Omega)}C_{H^2,L^\infty} 
+\|\kappa\|_{L^\infty(\Omega)}C_{H^1,L^4} \Bigr)
+\|\kappa\|_{L^\infty(\Omega)}C_{H^2,L^\infty}\Bigr) R < \frac12\,.
\end{equation}

Contractivity of $\mathcal{T}$ can be shown by taking $v^{(1)}$ and $v^{(2)}$ in $B_R$ and considering $u^{(1)}=\mathcal{T} v^{(1)}$ and $u^{(2)}=\mathcal{T} v^{(2)}$, whose differences $\overline{u}=u^{(1)} -u^{(2)} $ and $\overline{v}= v^{(1)} -v^{(2)} $ solve
\begin{equation} \label{eqn:contract}
(1-2\kappa v^{(1)})\overline{u}_{tt}+c^2 \mathcal{A} \overline{u}+b \mathcal{A}^\beta \partial_t^\alpha\overline{u}_t-2\kappa v_t^{(1)}\overline{u}_t=2\kappa\overline{v}_tu^{(2)}_t+2\kappa\overline{v}u^{(2)}_{tt}
\end{equation}
with homogeneous initial conditions. 
Similarly to above, with $\sigma=-2\kappa v^{(1)}$, $\mu=-2\kappa v^{(1)}_t$, $\rho(x,t)=0$, $h=2\kappa\overline{v}_tu^{(2)}_t+2\kappa\overline{v}u^{(2)}_{tt}$, since $v^{(1)}$, $u^{(2)}\in B_R$ (the latter due to the already shown self-mapping property of $\mathcal{T}$) we satisfy the conditions \eqref{eqn:nondegeneracy}, \eqref{eqn:nablasigmasmall}, \eqref{eqn:regularity_sigmamurhoh} on $\sigma$ and $\mu$. However $h$ in general fails to be contained in $L^2(0,T;\dot{H}^2(\Omega))$ (in particular the term $2\kappa\overline{v}u^{(2)}_{tt}$), hence we move to the lower order regularity regime from Proposition~\oldref{prop:wellposed_lin_CWK_low}.
To this end, we estimate 
\[
\|\nabla h\|_{\LtwoLtwo}
\leq
2\|\kappa\|_{L^4(\Omega)} \Bigl(
\|\overline{v}_t\|_{\LtwoLfourt} \|u^{(2)}_t\|_{\LinfLinft}
+\|\overline{v}\|_{\LinfLinft} \|u^{(2)}_{tt}\|_{\LtwoLfourt}\Bigr)\,.
\]
Thus imposing the additional smallness condition
\[
\theta := 2\|\kappa\|_{L^4(\Omega)} \sqrt{C} C_{H^1,L^4} C_{H^2,L^\infty}(\sqrt{T}+1) R < 1
\]
on $R$ and employing from Proposition~\oldref{prop:wellposed_lin_CWK_low}, we obtain contractivity
\[
\|\mathcal{T} v^{(1)}- \mathcal{T} v^{(2)}\|_{U_{lo}}=\|\overline{u}\|_{U_{lo}}\leq \theta \|\overline{v}\|_{U_{lo}}
= \theta \|v^{(1)}- v^{(2)}\|_{U_{lo}}
\]

\begin{theorem}\label{th:wellposedCKW}
For any $\alpha\in(0,1)$, $T>0$, $\kappa \in W^{2,4}(\Omega)\cap W^{1,\infty}(\Omega)$ there exists $R_0>0$ such that for any data $u_0\in\dot{H}^3(\Omega)$, $u_1\in\dot{H}^2(\Omega)$, $r\in L^2(0,T;\dot{H}^2(\Omega))$ satisfying
\begin{equation}\label{eqn:R0}
\|\mathcal{A}u_1\|_{\dot{L}^2(\Omega)}^2
+\|\mathcal{A}^{3/2} u_0\|_{\dot{L}^2(\Omega)}^2
+\|\mathcal{A}r\|_{\LtwoLtwod}^2
\leq R_0^2
\end{equation}
there exists a unique solution $u\in U$ of 
\begin{equation}\label{eqn:ivpCWK}
\begin{aligned}
(1-2\kappa u)u_{tt}+c^2\mathcal{A} u + b \mathcal{A}^\beta \partial_t^\alpha u &= 2\kappa (u_t)^2+r \mbox{ in }\Omega\times(0,T)\\
u(0)=u_0, \quad u_t(0)&=u_1 \mbox{ in }\Omega\,.
\end{aligned}
\end{equation}
\end{theorem}

\medskip

Existence of the linearisation of $G$ requires well-posedness of \eqref{eqn:Westervelt_init_D_lin} with $D = b \mathcal{A} \partial_t^\alpha$, that is,  \eqref{eqn:linPDE_CWK}, \eqref{eqn:init_CWK} with 
$\sigma=-2\kappa u$, $\mu=-4\kappa u_t$, $\rho=- 2\kappa u_{tt}$, $h=2\underline{\delta\kappa}(u \,u_{tt} + u_t^2)$. Due to the appearance of a $u_{tt}$ term we are in a similar situation to the contractivity proof above and therefore the lower regularity Proposition~\oldref{prop:wellposed_lin_CWK_low} is the right framework for analysing the linearisation of the forward problem.
\begin{proposition}
Under the assumptions of Theorem~\oldref{th:wellposedCKW}, for any $\underline{\delta\kappa}\in W^{1,\infty}(\Omega)$ there exists a unique solution $z\in U_{lo}$ of
\begin{equation}\label{eqn:ivpCWK_lin}
\begin{aligned}
(1-2\kappa u)z_{tt}+c^2\mathcal{A} z + b \mathcal{A}^\beta \partial_t^\alpha z - 4 \kappa u_t\, z_t - 2\kappa u_{tt} \, z
&= 2\underline{\delta\kappa}(u \,u_{tt} + u_t^2) \mbox{ in }\Omega\times(0,T)\\
z(0)=0, \quad z_t(0)&=0 \mbox{ in }\Omega\,,
\end{aligned}
\end{equation}
where $u\in U$ solves \eqref{eqn:ivpCWK}.
\end{proposition}

In order to prove Fr\'{e}chet differentiability we also need to bound the solution of \eqref{eqn:w}
that is,  \eqref{eqn:linPDE_CWK}, \eqref{eqn:init_CWK} with 
$\sigma=-2\kappa u$, $\mu=-4\kappa u_t$, $\rho=- 2\kappa u_{tt}$, 
$h=2\underline{d\kappa}(v \tilde{u}_{tt} + u v_{tt} + (\tilde{u}_t+u_t) v_t ) + 2\kappa (v v_{tt} + v_t^2)$, where $v$ can be bounded analogously to $z$ by Proposition~\oldref{prop:wellposed_lin_CWK_low}; in particular we can only expect to have $v_{tt}\in L^2(0,T; L^2(\Omega))$, so $h\in L^2(0,T;H^1(\Omega))$ is out of reach and we show Fr\'{e}chet differentiability in the very low regularity regime of Proposition~\oldref{prop:wellposed_lin_CWK_verylow}.

\begin{theorem}
For any $\alpha\in(0,1)$, $T>0$, $\bar{R}>0$ there exists $R_0>0$ such that for any data $u_0\in\dot{H}^3(\Omega)$, $u_1\in\dot{H}^2(\Omega)$, $r\in L^2(0,T;\dot{H}^2(\Omega))$ satisfying \eqref{eqn:R0}, the parameter-to-state map $G:B_{\bar{R}}(0)\to U$ is well-defined  according to Theorem \oldref{th:wellposedCKW}. Moreover, it is Fr\'{e}chet differentiable as an operator $G:B_{\bar{R}}(0)\to U_{vl}$.
Here $B_{\bar{R}}(0)=\{\kappa\in W^{2,4}(\Omega)\cap W^{1,\infty}(\Omega) \, : \, \|\kappa\|_{W^{2,4}(\Omega)\cap W^{1,\infty}(\Omega)}\leq\bar{R}\}$.
\end{theorem}

\subsection{Fractional Zener damping}\label{sec:analysisFZ}
Consider
\[
D = b_1 \mathcal{A} \partial_t^{\alpha_1}  + b_2 \partial_t^{\alpha_2+2} 
\quad \mbox{ with }b_2>0, \ b_1\geq b_2c^2, \  1\geq \alpha_1\geq\alpha_2>0,
\]
where based on the analysis in \cite{KaltenbacherRundell:2021b} we expect to get well-posedness of the nonlinear forward problem only in case $\alpha_1=1$, so we first of all focus on this case. Later on we will also prove a well-posedness result on the equation linearized at $\kappa=0$ in the practically relevant case $\alpha_1=\alpha_2=:\alpha$.

Again we first of all consider the initial boundary value problem for the
general linear {\sc pde}
\begin{equation}\label{eqn:linPDE_FZ}
(1-\sigma)u_{tt}+c^2\mathcal{A} u + b_1 \mathcal{A} u_t + b_2 \partial_t^{\alpha_2+2}u + \mu u_t + \rho u
= h
\end{equation}
\begin{equation}\label{eqn:init_FZ}
u(0)= u_0\,, \quad u_t(0)=u_1 \,, \quad (u_{tt}(0)=u_2 \mbox{ in case $\alpha_2>\frac12$)} 
\end{equation}
with given space- and time dependent functions $\sigma$, $\mu$, $\rho$, $h$.

Again we skip the details about the Faedo-Galerkin approach and the discretisation index $n$ and only provide the crucial energy estimate.
We multiply \eqref{eqn:linPDE_FZ} with $\mathcal{A}u_{tt}$ and integrate with respect to time, using the inequalities and identities
\[
\int_0^t\langle \mathcal{A} u(s), \mathcal{A}u_{tt}(s)\rangle\, ds
= \langle \mathcal{A} u(t), \mathcal{A}u_t(t)\rangle - \langle \mathcal{A} u(0), \mathcal{A}u_t(0)\rangle
-\int_0^t \|\mathcal{A}u_t(s)\|_{\dot{L}^2(\Omega)}^2
\]
and
\[
\begin{aligned}
&\int_0^t \langle \partial_t^{\alpha_2+2} [u](s), \mathcal{A}u_{tt}(s)\rangle\, ds
= \int_0^t \langle \partial_t^{\alpha_2} [\nabla u_{tt}](s), \nabla u_{tt}(s)\rangle\, ds
\\
&\geq\frac12 \int_0^t \partial_t^{\alpha_2}\left[\|\nabla u_{tt}\|_{\dot{L}^2(\Omega)}^2\right](s)\, ds
=\frac12 I_t^{1-\alpha}\left[\|\nabla u_{tt}\|_{L^2(\Omega)}^2\right](t)
\geq \frac{1}{2\Gamma(1-\alpha) t^\alpha} \|\nabla u_{tt}\|_{\LtwoLtwot}^2\,,
\end{aligned}
\]
where the latter equality holds provided $u_{tt}(0)=0$ 
and we have applied \eqref{eqn:Alikhanov_1} (actually to the Fourier components of the Galerkin discretisation $w=\lambda_j^{1/2} {u^n_j}''$) with $\gamma=\alpha_2$.

This yields the energy estimate
\[
\begin{aligned}
&\Bigl(\frac{b_2}{2\Gamma(1-\alpha) t^\alpha}+1-\overline{\sigma}\Bigr)
\|\nabla u_{tt}\|_{\LtwoLtwot}^2 
+ \frac{b_1}{2} \|\mathcal{A}u_t(t)\|_{\dot{L}^2(\Omega)}^2\\
&\leq \frac{b_1}{2} \|\mathcal{A}u_t(0)\|_{\LtwoLtwodt}^2
+ c^2 \langle \mathcal{A} u(0), \mathcal{A}u_t(0)\rangle\\ 
&\quad- c^2 \langle \mathcal{A} u(t), \mathcal{A}u_t(t)\rangle
+ c^2 \|\mathcal{A}u_t\|_{\LtwoLtwodt}^2 
+\int_0^t\langle \nabla\sigma u_{tt}+\nabla(\mu u_t + \rho u-h)(s),\nabla u_{tt}(s)\rangle\, ds\\
&\leq \frac{b_1}{2} \|\mathcal{A}u_t(0)\|_{\LtwoLtwodt}^2
+ c^2 \langle \mathcal{A} u(0), \mathcal{A}u_t(0)\rangle\\ 
&\qquad+ \frac{b_1}{4}\|\mathcal{A}u_t(t)\|_{\dot{L}^2(\Omega)}^2
+\frac{c^4}{b_1} \|\mathcal{A}u(t)\|_{\dot{L}^2(\Omega)}^2
+ c^2 \|\mathcal{A}u_t\|_{\LtwoLtwodt}^2\\
&\qquad+ 
\|\nabla \sigma\|_{\LinfLfourt} \|u_{tt}\|_{\LtwoLfourt}
\|\nabla u_{tt}\|_{\LtwoLtwot}
+ \frac{\epsilon}{2}\|\nabla u_{tt}\|_{\LtwoLtwot}^2 
\\
&\qquad+ \frac{1}{2\epsilon}\Bigl(
\|\nabla \mu\|_{\LinfLtwot} \|u_t\|_{L^2(0,t;L^\infty(\Omega)}
+\|\mu\|_{\LinfLfourt} \|\nabla u_t\|_{\LtwoLfourt}\\
&\qquad\qquad+
\|\nabla \rho\|_{\LtwoLtwot} \|u\|_{\LinfLinft}
+\|\rho\|_{\LtwoLfourt} \|\nabla u\|_{\LinfLfourt}
+ \|\nabla h\|_{\LtwoLtwot}\Bigr)^2
\end{aligned}
\]
Here we assume nondegeneracy 
\begin{equation}\label{eqn:nondegeneracy2}
\sigma(x,t)\leq \overline{\sigma}<\frac{b_2}{\Gamma(1-\alpha) T^\alpha}+1 \mbox{ for all }x\in\Omega\, \quad t\in(0,T)\,.
\end{equation}
and smallness of $\nabla\sigma$
\begin{equation}\label{eqn:nablasigmasmall2}
\|\nabla \sigma\|_{\LinfLfourt} < \frac{\frac{b_2}{\Gamma(1-\alpha) T^\alpha}+1-\overline{\sigma}}{C_{H^1\to L^4}}
\end{equation}
and choose $\epsilon< \frac{b_2}{\Gamma(1-\alpha) T^\alpha}+1\!-\!\underline{\sigma}-C_{H^1\!\to L^4}\|\nabla \sigma\|_{\LinfLfourt}$ to obtain,
using Gronwall's Lemma,
\begin{equation}\label{eqn:est_linFZ}
\begin{aligned}
&\|\nabla u_{tt}\|_{\LtwoLtwo}^2 
+ \|\mathcal{A}u_t\|_{\LinfLtwodt}^2
&\leq C \Bigl(\|\mathcal{A}u_t(0)\|_{\dot{L}^2(\Omega)}^2
+\|\nabla h\|_{\LtwoLtwodt}^2\Bigr)
\end{aligned}
\end{equation}
The required regularity on $\sigma$, $\mu$, $\rho$, $h$, $u_0$, $u_1$, is, besides \eqref{eqn:nondegeneracy2}, \eqref{eqn:nablasigmasmall2}
\begin{equation}\label{eqn:regularity_sigmamurhoh_FZ}
\begin{aligned}
&
\mu\in L^\infty(0,T;H^1(\Omega))\,, \quad \rho\in L^2(0,T;H^1(\Omega))\,, \\
&h\in L^2(0,T;\dot{H}^1(\Omega))\,, \quad
u_0, u_1\in \dot{H}^2(\Omega), \quad u_2=0\,.\end{aligned}
\end{equation}
\begin{proposition} \label{prop:wellposed_lin_FZ}
Under conditions \eqref{eqn:nondegeneracy2}, \eqref{eqn:nablasigmasmall2}, \eqref{eqn:regularity_sigmamurhoh_FZ}, there exists a unique solution 
\begin{equation}\label{defU_FZ}
u\in U:=H^2(0,T;\dot{H}^1(\Omega))\cap W^{1,\infty}(0,T;\dot{H}^2(\Omega))
\end{equation}
to the initial boundary value problem 
\eqref{eqn:linPDE_FZ}, \eqref{eqn:init_FZ}, and this solution satisfies the estimate \eqref{eqn:est_linFZ}.
\end{proposition}

\begin{theorem}\label{th:wellposedFZ}
For any $\alpha_2\in(0,1)$, $T>0$, $\kappa \in W^{1,4}(\Omega)$ there exists $R_0>0$ such that for any data $u_0,u_1\in\dot{H}^2(\Omega)$, $r\in L^2(0,T;\dot{H}^1(\Omega))$ satisfying
\begin{equation}\label{eqn:R0_FZ}
\|\mathcal{A}u_1\|_{\dot{L}^2(\Omega)}^2
+\|\mathcal{A} u_0\|_{\dot{L}^2(\Omega)}^2
+\|\nabla r\|_{\LtwoLtwo}^2
\leq R_0^2
\end{equation}
there exists a unique solution $u\in U$ of 
\begin{equation}\label{eqn:ivpFZ}
\begin{aligned}
b_2 \partial_t^{\alpha_2+2}u + (1-2\kappa u)u_{tt}+c^2\mathcal{A} u + b_1 \mathcal{A} u_t &= 2\kappa (u_t)^2+r \mbox{ in }\Omega\times(0,T)\\
u(0)=u_0, \quad u_t(0)=u_1, \quad u_{tt}&=0 \mbox{ in }\Omega\,.
\end{aligned}
\end{equation}
\end{theorem}

\begin{proposition}
Under the assumptions of Theorem~\oldref{th:wellposedFZ}, for any $\underline{\delta\kappa}\in L^3(\Omega)$ there exists a unique solution $z\in U$ of
\begin{equation}\label{eqn:ivpFZ_lin}
\begin{aligned}
b_2 \partial_t^{\alpha_2+2}z + (1\!-\!2\kappa u)z_{tt}\!+\!c^2\mathcal{A} z \!+\! b_1 \mathcal{A} z_t \!-\! 4 \kappa u_t\, z_t \!-\! 2\kappa u_{tt} \, z
&= 2\underline{\delta\kappa}(u \,u_{tt} + u_t^2) \mbox{ in }\Omega\times(0,T)\\
z(0)=0, \quad z_t(0), \quad z_{tt}=0&=0 \mbox{ in }\Omega\,.
\end{aligned}
\end{equation}
where $u\in U$ solves \eqref{eqn:ivpFZ}.
\end{proposition}

\begin{theorem}
For any $\alpha_2\in(0,1)$, $T>0$, $\bar{R}>0$ there exists $R_0>0$ such that for any data $u_0,u_1\in\dot{H}^2(\Omega)$, $r\in L^2(0,T;\dot{H}^2(\Omega))$ satisfying \eqref{eqn:R0_FZ}, the parameter-to-state map $G:B_{\bar{R}}(0)\to U$ is well-defined  according to Theorem \oldref{th:wellposedFZ}. Moreover, it is Fr\'{e}chet differentiable as an operator $G:B_{\bar{R}}(0)\to U$.
Here $B_{\bar{R}}(0)=\{\kappa\in W^{1,4}(\Omega) \, : \, \|\kappa\|_{W^{1,4}(\Omega)}\leq\bar{R}\}$.
\end{theorem}

We now consider the linear problem in case $\alpha_1=\alpha_2=:\alpha$, $b_1=b_2c^2+\delta$ with $\delta\geq0$
\begin{equation}\label{eqn:linPDE_FZ_samealpha}
(1-\sigma)u_{tt}+c^2\mathcal{A} u + b_1 \mathcal{A} \partial_t^\alpha u + b_2 \partial_t^{\alpha+2}u + \mu u_t + \rho u
= h
\end{equation}
in which the differential operator can partially be factorised as 
\[
\begin{aligned}
&(1-\sigma)\partial_{tt}+c^2\mathcal{A}  + (b_2c^2+\delta) \mathcal{A} \partial_t^\alpha  + b_2 \partial_t^{\alpha+2} + \mu \partial_t + \rho \mbox{id}\\
&= \Bigl(\partial_{tt}+c^2\mathcal{A}\Bigr)\Bigl(b_2\partial_t^\alpha  + \mbox{id}\Bigr)
-\sigma \partial_{tt} + \delta \mathcal{A} \partial_t^\alpha + \mu \partial_t + \rho \mbox{id}
\end{aligned}
\]
Thus, up to the ``perturbation'' terms containing $\sigma$, $\delta$, $\mu$, and $\rho$, the auxiliary function $\tilde{u}=b_2\partial_t^\alpha u + u$ satisfies a wave equation $\tilde{u}_{tt}+c^2\mathcal{A}\tilde{u}=h$. Motivated by this fact, we multiply \eqref{eqn:linPDE_FZ_samealpha} with $\mathcal{A}\tilde{u}_t$ to obtain the energy identity 
\begin{equation}\label{eqn:enid_fZ_lin}
\begin{aligned}
\frac12\|\nabla (b_2\partial_t^\alpha u + u)_t(t)\|_{L^2(\Omega)}^2& 
+\frac{c^2}{2}\|\mathcal{A} (b_2\partial_t^\alpha u + u)(t)\|_{\dot{L}^2(\Omega)}^2\\
&\qquad+ \delta\int_0^t \langle \nabla \partial_t^\alpha u(s) , \nabla(b_2\partial_t^\alpha u + u)_t(s)\rangle\, ds\\ 
= \frac12\|\nabla (b_2\partial_t^\alpha u + u)_t(0)&\|_{L^2(\Omega)}^2 
+\frac{c^2}{2}\|\mathcal{A} (b_2\partial_t^\alpha u + u)(0)\|_{\dot{L}^2(\Omega)}^2\\
&\quad+\int_0^t \langle \nabla(h+\sigma u_tt-\mu u_t-\rho u)(s), \nabla(b_2\partial_t^\alpha u + u)_t(s)\rangle\, ds\,.
\end{aligned}
\end{equation}
The term containing $\delta$ can be nicely tackled by means of Lemma \oldref{lem:Alikhanov2}
\[
\begin{aligned}
\int_0^t \langle \nabla \partial_t^\alpha u(s)& , \nabla(b_2\partial_t^\alpha u + u)_t(s)\rangle\, ds\\
&= \frac{b_2}{2} \|\nabla \partial_t^\alpha u(t)\|_{L^2(\Omega)}^2 - \frac{b_2}{2} \|\nabla \partial_t^\alpha u(0)\|_{L^2(\Omega)}^2
+ \int_0^t \langle \nabla \partial_t^\alpha u(s) , \nabla u_t(s)\rangle\, ds\\
&\geq \frac{b_2}{2} \|\nabla \partial_t^\alpha u(t)\|_{L^2(\Omega)}^2
+\frac{1}{2\Gamma(\alpha)t^{1-\alpha}}\|\nabla\partial_t^{\alpha} u\|_{\LtwoLtwot}^2\,,
\end{aligned}
\]
(which reflects the physical fact that $\delta$ is the diffusivity of sound and therefore the corresponding term models damping).
However, in the term containing $\sigma$ this is inhibited by the time-dependence of $\sigma$. 
Thus in case $\alpha_1=\alpha_2=:\alpha<1$ we have to restrict ourselves to the linearisation of the forward problem at $\kappa=0$ (where also $\mu=0$, $\rho=0$), where the above together with Young's inequality yields the energy estimate
\begin{equation}\label{eqn:enest_fZ_lin}
\begin{aligned}
\frac12\|\nabla &(b_2\partial_t^\alpha u + u)_t\|_{\LinfLtwot}^2 
+c^2\|\mathcal{A} (b_2\partial_t^\alpha u + u)\|_{\LinfLtwodt}^2\\
&+ \delta b_2 \|\nabla \partial_t^\alpha u(t)\|_{\LinfLtwot}^2
+\frac{1}{2\Gamma(\alpha)t^{1-\alpha}}\|\nabla\partial_t^{\alpha} u\|_{\LtwoLtwot}^2
\\ 
&\leq \|\nabla (b_2\partial_t^\alpha u + u)_t(0)\|_{L^2(\Omega)}^2 
+c^2\|\mathcal{A} (b_2\partial_t^\alpha u + u)(0)\|_{\dot{L}^2(\Omega)}^2
+ \frac12\|\nabla h\|_{\LoneLtwot}^2 
\end{aligned}
\end{equation}
The {\sc pde} yields an estimate of $\tilde{u}_{tt}$ as follows
\begin{equation}\label{eqn:enest_fZ_lin_tt}
\begin{aligned}
&\|(b_2\partial_t^\alpha u + u)_{tt}\|_{\LinfLtwot} 
=\|h-c^2\mathcal{A} (b_2\partial_t^\alpha u + u)\|_{\LinfLtwodt}^2\\
\end{aligned}
\end{equation}
To extract temporal reguarity of $u$ from regularity of $b_2\partial_t^\alpha u + u$ for $b_2>0$, we make use of known regularity results of time fractional ODEs: $b_2\partial_t^\alpha u + u =\tilde{u}\in W^{k,\infty}(0,T;Z)$ implies $u\in W^{k+\alpha,\infty}(0,T;Z)$.

In order to formulate a result on this linearisation, we will postulate the regularity of $u$ that is to be  ``expected'' from \eqref{eqn:enest_fZ_lin}, \eqref{eqn:enest_fZ_lin_tt} (but apparently not provable): 
\begin{equation}\label{defU_FZ_samealpha}
u\in U:=W^{2+\alpha,\infty}(0,T;\dot{L}^2(\Omega))\cap W^{1+\alpha,\infty}(0,T;\dot{H}^1(\Omega))\cap W^{\alpha,\infty}(0,T;\dot{H}^2(\Omega))
\end{equation}
The required regularity of $\underline{\delta\kappa}$ to guarantee $h=2\underline{\delta\kappa}(u \,u_{tt} + u_t^2)\in L^1(0,t;L^2(\Omega))$ is therefore $\underline{\delta\kappa}\in L^\infty(\Omega)$.

We thus obtain the following result on the linearisation of the forward problem at $\kappa=0$:

\begin{proposition}\label{prop:FZalpha12}
For any $u\in U$ and $\underline{\delta\kappa}\in L^\infty(\Omega)$ there exists a unique solution $z\in U$ of
$$
\begin{aligned}
b_2 \partial_t^{\alpha+2}z + z_{tt}+c^2\mathcal{A} z + b_1 \mathcal{A} \partial_t^\alpha z
&= 2\underline{\delta\kappa}(u \,u_{tt} + u_t^2) \mbox{ in }\Omega\times(0,T)\\
z(0)=0, \quad z_t(0), \quad z_{tt}=0&=0 \mbox{ in }\Omega\,.
\end{aligned}
$$
\end{proposition}

\section{Injectivity of the linearised forward operator}\label{sect:inject}

The forwards map is defined by $F(\kappa)=\mbox{tr}_\Sigma u$, where $\mbox{tr}_\Sigma v$ denotes the time trace of the space and time dependent function $v:(0,T)\times \Omega$ at the observation surface $\Sigma$ (which may also just be a single point $\Sigma=\{x_0\}$) and $u$ solves
\begin{equation}\label{eqn:Westervelt_init_D_uni}
\begin{aligned}
u_{tt}+c^2\mathcal{A} u + D u &= \kappa(x)(u^2)_{tt} + r \mbox{ in }\Omega\times(0,T)\\
u(0)=0, \quad u_t(0)&=0 \mbox{ in }\Omega\,.
\end{aligned}
\end{equation}
Its linearisation at $\kappa=0$ in direction $\underline{\delta\kappa}$ is $F'(0)\underline{\delta\kappa}=\mbox{tr}_\Sigma z_0$, where $z_0$ solves
\begin{equation}\label{eqn:z0}
z_{tt}+c^2\mathcal{A} z + D z = \underline{\delta\kappa}(u_0^2)_{tt} 
\end{equation}
where 
\begin{equation}\label{eqn:u0}
u_{0,tt}+c^2\mathcal{A} u_0 + D u_0 = r. 
\end{equation}
Both {\sc pde}s \eqref{eqn:z0}, \eqref{eqn:u0} come with homogeneous initial conditions.

As in the previous section, we consider the two damping models
\begin{equation}\label{eqn:CHCW_uni}
D = b \mathcal{A}^\beta \partial_t^\alpha 
\quad \mbox{ ({\sc cwch})}
\end{equation}
and
\begin{equation}\label{eqn:FZ_uni}
D = b_1 \mathcal{A} \partial_t^{\alpha_1}  + b_2 \partial_t^{\alpha_2+2} 
\quad \mbox{ ({\sc fz})}\,.
\end{equation}

The Laplace transformed solutions to the corresponding resolvent equation
\[
\hat w(\lambda,s) = \frac{1}{\omega(\lambda,s)} \mbox{ with }
\omega(\lambda,s)= \begin{cases} 
s^2 + b\lambda^\beta s^\alpha + c^2\lambda \mbox{ for {\sc cwch}}\\
b_2 s^{2+\alpha_2} + s^2 + b_1\lambda s^{\alpha_1} + c^2\lambda \mbox{ for {\sc fz}}
\end{cases}
\]
will play a crucial role in the proofs below.

From \cite[Lemma 4.1]{KaltenbacherRundell:2021b} we conclude that in case of
{\sc cwch}, the function $\omega^{{\sc cwch}}(\lambda,\cdot)$ has precisely two complex-conjugate zeros $p^{{\sc cwch}}_{+}(\lambda)$, $p^{{\sc cwch}}_{-}(\lambda)$, which lie in the left hand complex plane.

For {\sc fz}, we first consider the particular parameter configuration (corresponding to vanishing viscosity of sound) 
\begin{equation}\label{eqn:b12alpha12}
b_1 = b_2c^2 \mbox{ and }\alpha_1=\alpha_2
\end{equation}
in which we can factorise 
$\omega^{FZ}(\lambda,s) = (b_2 s^{\alpha_2}+1)(s^2+c^2\lambda)$ and get the roots
\[
p^{FZ}_0=-\frac{1}{b_2}\mbox{ (only in case $\alpha_1=\alpha_2=1$)}\,, \qquad
p^{FZ}_\pm(\lambda)=\pm i c\sqrt{\lambda}
\]
Note that $p^{FZ}_0$ is independent of $\lambda$, but $p^{FZ}_\pm(\lambda)$
obviously allows to distinguish between different $\lambda$'s.
This distinction is possible in general, a fact that has already been shown for the {\sc cwch} case with $\beta=1$ in \cite[Remark 4.1]{KaltenbacherRundell:2021b}.
As an additional result, that is not needed for the uniqueness proof but might be convenient for the computation of poles and residues, we state that the poles are single in certain cases.
\begin{lemma}\label{lem:what}
The poles of $\hat{w}^{{\sc cwch}}$ and of $\hat{w}^{FZ}$ (except for $p^{FZ}_0$ in case $\alpha_1=\alpha_2=1$) differ for different $\lambda$.  
Moreover, in the case {\sc cwch} and in the case {\sc fz} with \eqref{eqn:b12alpha12} the poles are single.
\end{lemma}
\begin{proof}
For {\sc cwch}, let $f(z) = z^2 + c^2\lambda$, $g(z) = b\lambda^{\beta} z^\alpha$.
Then for a sufficiently large $R>c^2\lambda$ let
$C_R$ be the circle radius $R$, centre at the origin.
Then $|g(z)| < |f(z)|$ on $C_R$ and so Rouch\'e's theorem shows that
$f(z)$ and $(f+g)(z)$ have the same number of roots, counted with multiplicity,
within $C_R$.
For $f$ these are only at $z=\pm i\sqrt{\lambda}c$
so the same must be true of $f+g$
and so $\omega^{{\sc cwch}}$ has 
precisely one single root in the third and in the fourth quadrant, respectively.

Suppose now that $\hat{w}^{{\sc cwch}}$ has a pole at $r e^{i\theta}$,
where $\pi/2<\theta<\pi$,
for both $\lambda_1$ and $\lambda_2$. Then for $s=r e^{i\theta}$
$$
s^2 + b\lambda_1^{\beta} s^\alpha + c^2\lambda_1 = 0 \qquad
s^2 + b\lambda_2^{\beta} s^\alpha + c^2\lambda_2 = 0
$$
so that 
$$
\frac{b}{c^2}
\frac{(\lambda_1^{\beta}-\lambda_2^{\beta})}{(\lambda_1-\lambda_2)} = -s^\alpha
$$
Now if $\lambda_1\not=\lambda_2$ then the left hand side is positive and real
and so $\alpha\theta =\pi$.
This means that $\theta > \pi$, a contradiction.

In case of {\sc fz}, assuming that $p$ is a pole of both $\hat{w}^{FZ}(\lambda_1,\cdot)$ and $\hat{w}^{FZ}(\lambda_2,\cdot)$
we have 
\[
0 = \omega(\lambda_1,p)-\omega(\lambda_2,p)
=(\lambda_1-\lambda_2)[b_1 p^{\alpha_1}+c^2]\,,
\]
where due $\omega(\lambda_1,p)=0$, the term in brackets $\,b_1 p^{\alpha_1}+c^2=-\frac{p^2}{\lambda_1}(b_2p^{\alpha_2}+1)\not=0\,$, hence $\lambda_1=\lambda_2$.
In the factorisable case \eqref{eqn:b12alpha12} of {\sc fz},
obviously all roots are single.
\end{proof}

As in \cite{KaltenbacherRundell:2021c} (where we used the classical damping term $D=b\mathcal{A}\partial_t$), we assume that $r$
has the form 
\begin{equation}\label{eqn:choicer}
r(x,t) = f(x) \chi''(t) + c^2\mathcal{A}f(x)\chi(t) + D[f(x)\chi(t)]
\end{equation}
with some function $f$ in the domain of $\mathcal{A}$ vanishing only on a set of 
measure zero and some twice differentiable function $\chi$ of time such that $(\chi^2)''(t_0)\not=0$ for some $t_0>0$.
With \eqref{eqn:choicer}, the solution $u_0$ of equation~\eqref{eqn:u0} is clearly given by 
$u_0(x,t) = f(x) \chi(t)$, so that
$\underline{d\kappa}(u_0^2)_{tt}$ can be written in the form  
\begin{equation}\label{eqn:u0_fac2}
\bigl(\underline{d\kappa}(x)u_0^2(x,t)\bigr)_{tt} = \sum_{j=1}^\infty a_j\phi_j(x)\psi_j(t)  
\end{equation}
where $a_j$ are the coefficients of $\underline{d\kappa}\cdot f$ with respect to the eigenfunction basis $(\phi_j)_{j\in\mathbb{N}}$, and $\psi_j=(\chi^2)''$.

We can rewrite equation~\eqref{eqn:z0} as
\begin{equation}\label{eqn:z0j2}
z_j''(t)+c^2\lambda_j z_j(t) + D_j z_j
 = a_j \psi_j(t)\,, \quad t>0\,, \qquad z_j(0)=0\,, \ z_j'(0)=0
\end{equation}
for all $j\in\mathbb{N}$, where
\[
z_0(x,t)=\sum_{j=1}^\infty z_j(t)\phi_j(x)\,, \qquad 
D_j = \begin{cases}b \lambda_j^\beta \partial_t^\alpha  \mbox{ for {\sc cwch}}\\
b_2 \partial_t^{2+\alpha_2} + b_1\lambda_j \partial_t^{\alpha_1} \mbox{ for {\sc fz}}
\end{cases}
\]
Applying the Laplace transform to both sides of \eqref{eqn:z0j2} yields
\begin{equation}\label{eqn:z0jhat2}
\hat{z}_j(s)= \hat w_j(s) a_j \widehat{\psi}_j(s), \qquad
\mbox{where }\ \hat w_j(s) = \frac{1}{\omega(\lambda_j,s)}
\,, \qquad s\in\mathbb{C}\,,
\end{equation}
and we have used homogeneity of the initial conditions.

Thus, assuming that $F'(0)\underline{d\kappa} = \mbox{tr}_\Sigma z_0=0$ implies that 
\[
0 = \hat{z}_0(x_0,s) = \sum_{j=1}^\infty a_j \phi_j(x_0) \hat{w}_j(s) \widehat{\psi}_j(s)\,, \qquad 
\mbox{ for all }s\in\mathbb{C}\,, \ x_0\in  \Sigma\,.
\]
Considering the residues at some pole $p_m$ corresponding to the eigenvalue $\lambda_m$ and using the fact that by Lemma~\oldref{lem:what}, $(s-p_m) \hat{w}(\lambda_j,s)=0$ for $j\not=m$ yields
\[
\begin{aligned}
0 &= \mbox{Res}(\hat{z}_0(x_0;p_m))
= \sum_{j=1}^\infty a_j \phi_j(x_0) \lim_{s\to p_m}(s-p_m)^{\ell_m} \hat{w}(\lambda_j,s) \widehat{\psi}_j(s)\\ 
&= \mbox{Res}(\hat{w}_m;p_m) \widehat{\psi}_j(p_m) \sum_{k\in K_m} a_k \phi_k(x_0) 
\,.
\end{aligned}
\]
Here $\ell_m$ is the multiplicity of $p_m$ as a root of $\omega(\lambda_m,\cdot)$ and $K_m\subseteq\mathbb{N}$ is an enumeration of the eigenspace basis $(\phi_k)_{k\in K_m}$ corresponding to the eigenvalue $\lambda_m$.
Assuming now that
\begin{equation}\label{eqn:ass_inj_psi}
\widehat{\psi}_k(p_m)\not=0 
\end{equation}
and there exists points $x_{0,m,1}, \ldots x_{0,m,N_m}\in \Sigma$, $N_m\geq \# K_m$ such that 
\begin{equation}\label{eqn:ass_inj_Sigma}
\mbox{the matrix }\phi_k(x_{0,m,i})_{k\in K_m,i\in\{1,\ldots,N_m\}} \mbox{ has full rank }\# K_m
\end{equation}
we can conclude that $a_k=0$ for all $k\in K_m$.

Now since $(u_0^2)_{tt}(t_0)=f (\chi^2)''(t_0)$ only vanishes on a set of measure zero and \eqref{eqn:u0_fac2}, we can conclude that $\underline{d\kappa}=0$ almost everywhere. 

\begin{theorem}\label{th:lininj}
Under the above assumptions \eqref{eqn:u0_fac2}, \eqref{eqn:ass_inj_psi}, \eqref{eqn:ass_inj_Sigma} for all $m\in\mathbb{N}$, $k\in K_m$, the linearised derivative at $\kappa=0$, $F'(0)$ is injective.
\end{theorem}

In particular, \eqref{eqn:ass_inj_Sigma} is satisfied in the spatially 1-dimensional case $\Sigma=\{x_0\}$, where all eigenvalues of $\mathcal{A}$ are single, i.e., $\# K_m=1$ for all $m$, provided none of the eigenfunctions vanish at $x_0$; this can be achieved by taking $x_0$ on the boundary and where $\phi_j$ is subject to non-Dirichlet conditions.

\section{Ill-posedness of the linearised inverse problem}

As in the injectivity section, we consider the linear(ised at $\kappa=0$) problem of recovering $\underline{\delta\kappa}(x)$ from time trace observations 
$$
h(t) = \mbox{tr}_\Sigma z_0 = z_0(x_0,t_0), \quad x_0\in\Sigma, 
$$
where $z_0$ solves 
\begin{equation}\label{eqn:z0_illp}
z_{tt}+c^2\mathcal{A} z + D z = \underline{\delta\kappa}(u_0^2)_{tt} 
\end{equation}
with $u_0$ solving
\begin{equation}\label{eqn:u0_illp}
u_{0,tt}+c^2\mathcal{A} u_0 + D u_0 = r 
\end{equation}
both with homogeneous initial conditions.

Again, we assume that the excitation $r$ has been chosen such that $u_0$ takes the form $u_0(x,t) = f(x) \chi(t)$ and employ the shorthand notation $\tilde{\chi}=(\chi^2)''$.
Using the eigensystem of $\mathcal{A}$ we can then write
$$
z_0(x,t) = \sum_{j=1}^\infty \phi(x) z_j(t)
$$
where
\begin{equation}\label{eqn:z0jhat3}
\hat{z}_j(s)= \hat w_j(s) \langle\underline{\delta\kappa}\cdot f,\phi_j\rangle  \widehat{\tilde{\chi}}(s), \qquad \hat w_j(s) = \frac{1}{\omega(\lambda_j,s)}
\,, \qquad s\in\mathbb{C}\,,
\end{equation}
and 
\[
\omega(\lambda,s)= \begin{cases} 
s^2 + b\lambda^\beta s^\alpha + c^2\lambda \mbox{ for {\sc cwch}}\\
b_2 s^{2+\alpha_2} + s^2 + b_1\lambda s^{\alpha_1} + c^2\lambda \mbox{ for {\sc fz}}\,.
\end{cases}
\]
As in the injectivity section we obtain (for simplicity in the 1-d case where all eigenvalues are single)
\[
\begin{aligned}
\mbox{Res}(\hat{h}(p_m))
= \mbox{Res}(\hat{w}_m;p_m) \widehat{\tilde{\chi}}(p_m) \langle\underline{\delta\kappa}\cdot f,\phi_m\rangle \phi_m(x_0) 
\,.
\end{aligned}
\]
that is, 
\begin{equation}\label{eqn:deltakappafphim}
\begin{aligned}
\langle\underline{\delta\kappa}\cdot f,\phi_m\rangle
=\mbox{Res}(\hat{h}(p_m))
\left(\mbox{Res}(\hat{w}_m;p_m) \widehat{\tilde{\chi}}(p_m)  \phi_m(x_0) \right)^{-1}
\,.
\end{aligned}
\end{equation}
(In our numerical example we had $\chi(t)=t$, hence $\widehat{\tilde{\chi}}(p_m)=2$.)

By l'Hospital's rule we have  
\[
\mbox{Res}(\hat{w}_m;p_m)= \lim_{s\to p_m} \frac{s-p_{m}}{\omega(\lambda,s)} =
\lim_{s\to p_m} \frac{1}{\omega'(\lambda,s)} =
\begin{cases} 
\frac{1}{2p_m + \alpha b\lambda_j^\beta p_m^{\alpha-1}} \mbox{ for {\sc cwch}}\\
\frac{1}{(2+\alpha_2)b_2 p_m^{1+\alpha_2} + 2 p_m + \alpha_1b_1\lambda_j p_m^{\alpha_1-1}} \mbox{ for {\sc fz}}\,.
\end{cases}
\]
Thus, the factor multiplied with $\langle\underline{\delta\kappa}\cdot f,\phi_m\rangle$ in \eqref{eqn:deltakappafphim} only mildly grows with $p_m$.

The major ill-posedness seems to lie in the evaluation of the residue of the observations $\mbox{Res}(\hat{h}(p_m))$ at the poles $p_m$, form knowledge of $h(t)$ for $t>0$, that is, from $\hat{h}(s)=\int_0^\infty e^{-st} h(t)\, dt$ for $s$ with nonnegative real part (so that the integral defining the Laplace transform is well-defined).  If these poles lie on the imaginary axis (wave equation), this is still well posed. The further left the poles lie, the more ill-posed this problem.

\font\smallsymbol = cmmi8
\newdimen\xfiglen \newdimen\yfiglen
\xfiglen=2.5 true in
\yfiglen=10 true in
%
%
\setbox\figurelegendone=\hbox{
\beginpicture
  \setcoordinatesystem units <0.08\xfiglen,0.012\yfiglen> 
  \setplotarea x from 0 to 2, y from 0 to 4
\linethickness=0.6pt
\scriptsize
   \setsolid  
  \put {$\star$\ $\alpha=0.25$}  [l] at 0 3
  \put {$\diamond$\ $\alpha=0.5$}  [l] at 0 2
  \put {$\circ$\ $\alpha=0.9$}  [l] at 0 1
  \put {$\bullet$\ $\alpha=1.0$} [l] at 0 0
\endpicture
}
\setbox\figurelegendtwo=\hbox{
\beginpicture
  \setcoordinatesystem units <0.07\xfiglen,0.014\yfiglen> 
  \setplotarea x from 0 to 2, y from 0 to 2
\linethickness=0.6pt
\scriptsize
   \setsolid  
  \put {$\star$\ $\;\beta=1$}  [l] at 0 2
  \put {$\diamond$\ $\;\beta=\frac{3}{4}$}  [l] at 0 1
  \put {$\circ$\ $\;\beta=\frac{1}{2}$}  [l] at 0 0
\endpicture
}
\setbox\figurelegendthree=\hbox{
\beginpicture
  \setcoordinatesystem units <0.1\xfiglen,0.012\yfiglen> 
  \setplotarea x from 0 to 2, y from 0 to 3
\linethickness=0.6pt
\scriptsize
   \setsolid  
  \put {$\star$\ $\;\delta=0.1$}  [l] at 0 3
  \put {$\diamond$\ $\;\delta=0.5$}  [l] at 0 2
  \put {$\circ$\ $\;\delta=1.0$}  [l] at 0 1
  \put {$\bullet$\ $\;\delta=2.0$}  [l] at 0 0
\endpicture
}
\setbox\figureone=\vbox{\hsize=\xfiglen
\beginpicture
  \setcoordinatesystem units <\xfiglen,\yfiglen>  point at 0.0 0.0
  \setplotarea x from 0 to 1, y from 0 to 0.17
\scriptsize
  \axis bottom shiftedto y=0.0 ticks short numbered from 0 to 1 by 0.2 /
  \axis left ticks short numbered from 0 to 0.17 by 0.05 /
\footnotesize
\put {$\kappa(x)$} [lt] at 0.02 0.17
\put {$\alpha=1.0$} [lt] at 0.4 0.16
\setlinear
\setsolid
\Red{\relax   
%
\plot
         0         0
    0.0100         0
    0.0200         0
    0.0300         0
    0.0400         0
    0.0500         0
    0.0600    0.0100
    0.0700    0.0200
    0.0800    0.0300
    0.0900    0.0400
    0.1000    0.0500
    0.1100    0.0600
    0.1200    0.0700
    0.1300    0.0800
    0.1400    0.0900
    0.1500    0.1000
    0.1600    0.1060
    0.1700    0.1120
    0.1800    0.1180
    0.1900    0.1240
    0.2000    0.1300
    0.2100    0.1300
    0.2200    0.1300
    0.2300    0.1300
    0.2400    0.1300
    0.2500    0.1300
    0.2600    0.1200
    0.2700    0.1100
    0.2800    0.1000
    0.2900    0.0900
    0.3000    0.0800
    0.3100    0.0750
    0.3200    0.0700
    0.3300    0.0650
    0.3400    0.0600
    0.3500    0.0550
    0.3600    0.0500
    0.3700    0.0450
    0.3800    0.0400
    0.3900    0.0350
    0.4000    0.0300
    0.4100    0.0315
    0.4200    0.0330
    0.4300    0.0345
    0.4400    0.0360
    0.4500    0.0375
    0.4600    0.0390
    0.4700    0.0405
    0.4800    0.0420
    0.4900    0.0435
    0.5000    0.0450
    0.5100    0.0480
    0.5200    0.0510
    0.5300    0.0540
    0.5400    0.0570
    0.5500    0.0600
    0.5600    0.0630
    0.5700    0.0660
    0.5800    0.0690
    0.5900    0.0720
    0.6000    0.0750
    0.6100    0.0825
    0.6200    0.0900
    0.6300    0.0975
    0.6400    0.1050
    0.6500    0.1125
    0.6600    0.1200
    0.6700    0.1275
    0.6800    0.1350
    0.6900    0.1425
    0.7000    0.1500
    0.7100    0.1500
    0.7200    0.1500
    0.7300    0.1500
    0.7400    0.1500
    0.7500    0.1500
    0.7600    0.1500
    0.7700    0.1500
    0.7800    0.1500
    0.7900    0.1500
    0.8000    0.1500
    0.8100    0.1450
    0.8200    0.1400
    0.8300    0.1350
    0.8400    0.1300
    0.8500    0.1250
    0.8600    0.1200
    0.8700    0.1150
    0.8800    0.1100
    0.8900    0.1050
    0.9000    0.1000
    0.9100    0.0800
    0.9200    0.0600
    0.9300    0.0400
    0.9400    0.0200
    0.9500         0
    0.9600         0
    0.9700         0
    0.9800         0
    0.9900         0
    1.0000         0
/\relax

}\relax
\setdashes <2pt>
\Black{\relax 

}\relax
\setplotsymbol ({\sixrm .})
\setsolid
\Blue{\relax   
\plot
         0    0.0000
    0.0100    0.0002
    0.0200    0.0005
    0.0300    0.0011
    0.0400    0.0022
    0.0500    0.0033
    0.0600    0.0059
    0.0700    0.0086
    0.0800    0.0122
    0.0900    0.0167
    0.1000    0.0213
    0.1100    0.0277
    0.1200    0.0341
    0.1300    0.0412
    0.1400    0.0490
    0.1500    0.0568
    0.1600    0.0650
    0.1700    0.0733
    0.1800    0.0811
    0.1900    0.0886
    0.2000    0.0960
    0.2100    0.1015
    0.2200    0.1069
    0.2300    0.1108
    0.2400    0.1132
    0.2500    0.1156
    0.2600    0.1144
    0.2700    0.1132
    0.2800    0.1102
    0.2900    0.1055
    0.3000    0.1008
    0.3100    0.0934
    0.3200    0.0860
    0.3300    0.0780
    0.3400    0.0694
    0.3500    0.0607
    0.3600    0.0524
    0.3700    0.0440
    0.3800    0.0367
    0.3900    0.0302
    0.4000    0.0238
    0.4100    0.0205
    0.4200    0.0171
    0.4300    0.0155
    0.4400    0.0157
    0.4500    0.0158
    0.4600    0.0192
    0.4700    0.0225
    0.4800    0.0270
    0.4900    0.0326
    0.5000    0.0382
    0.5100    0.0447
    0.5200    0.0512
    0.5300    0.0576
    0.5400    0.0638
    0.5500    0.0700
    0.5600    0.0752
    0.5700    0.0804
    0.5800    0.0851
    0.5900    0.0892
    0.6000    0.0933
    0.6100    0.0969
    0.6200    0.1004
    0.6300    0.1041
    0.6400    0.1077
    0.6500    0.1114
    0.6600    0.1157
    0.6700    0.1200
    0.6800    0.1245
    0.6900    0.1293
    0.7000    0.1340
    0.7100    0.1385
    0.7200    0.1430
    0.7300    0.1468
    0.7400    0.1499
    0.7500    0.1531
    0.7600    0.1541
    0.7700    0.1551
    0.7800    0.1550
    0.7900    0.1536
    0.8000    0.1522
    0.8100    0.1487
    0.8200    0.1451
    0.8300    0.1405
    0.8400    0.1349
    0.8500    0.1293
    0.8600    0.1214
    0.8700    0.1135
    0.8800    0.1042
    0.8900    0.0935
    0.9000    0.0827
    0.9100    0.0695
    0.9200    0.0563
    0.9300    0.0431
    0.9400    0.0298
    0.9500    0.0165
    0.9600    0.0067
    0.9700         0
    0.9800         0
    0.9900    0.0003
    1.0000    0.0011
/\relax

}\relax
\endpicture
}
\setbox\figuretwo=\vbox{\hsize=\xfiglen
\beginpicture
  \setcoordinatesystem units <\xfiglen,\yfiglen>  point at 0.0 0.0
  \setplotarea x from 0 to 1, y from 0 to 0.17
\scriptsize
  \axis bottom shiftedto y=0.0 ticks short numbered from 0 to 1 by 0.2 /
  \axis left ticks short numbered from 0 to 0.17 by 0.05 /
\footnotesize
\put {$\kappa(x)$} [lt] at 0.02 0.17
\put {$\alpha=0.9$} [lt] at 0.4 0.16
\setlinear
\setsolid
\Red{\relax   

}\relax
\setdashes <2pt>
\Black{\relax 

}\relax
\setplotsymbol ({\sixrm .})
\setsolid
\Blue{\relax   
\plot
         0    0.0000
    0.0100    0.0003
    0.0200    0.0007
    0.0300    0.0016
    0.0400    0.0031
    0.0500    0.0046
    0.0600    0.0083
    0.0700    0.0119
    0.0800    0.0167
    0.0900    0.0228
    0.1000    0.0289
    0.1100    0.0370
    0.1200    0.0452
    0.1300    0.0539
    0.1400    0.0631
    0.1500    0.0724
    0.1600    0.0812
    0.1700    0.0901
    0.1800    0.0980
    0.1900    0.1048
    0.2000    0.1117
    0.2100    0.1151
    0.2200    0.1186
    0.2300    0.1200
    0.2400    0.1194
    0.2500    0.1187
    0.2600    0.1140
    0.2700    0.1093
    0.2800    0.1031
    0.2900    0.0953
    0.3000    0.0875
    0.3100    0.0781
    0.3200    0.0688
    0.3300    0.0596
    0.3400    0.0505
    0.3500    0.0414
    0.3600    0.0344
    0.3700    0.0273
    0.3800    0.0218
    0.3900    0.0180
    0.4000    0.0141
    0.4100    0.0139
    0.4200    0.0138
    0.4300    0.0153
    0.4400    0.0185
    0.4500    0.0217
    0.4600    0.0270
    0.4700    0.0324
    0.4800    0.0381
    0.4900    0.0442
    0.5000    0.0503
    0.5100    0.0556
    0.5200    0.0609
    0.5300    0.0655
    0.5400    0.0693
    0.5500    0.0731
    0.5600    0.0756
    0.5700    0.0782
    0.5800    0.0805
    0.5900    0.0828
    0.6000    0.0851
    0.6100    0.0884
    0.6200    0.0917
    0.6300    0.0959
    0.6400    0.1011
    0.6500    0.1063
    0.6600    0.1130
    0.6700    0.1197
    0.6800    0.1265
    0.6900    0.1333
    0.7000    0.1401
    0.7100    0.1452
    0.7200    0.1502
    0.7300    0.1538
    0.7400    0.1558
    0.7500    0.1578
    0.7600    0.1567
    0.7700    0.1556
    0.7800    0.1535
    0.7900    0.1505
    0.8000    0.1474
    0.8100    0.1437
    0.8200    0.1400
    0.8300    0.1361
    0.8400    0.1319
    0.8500    0.1277
    0.8600    0.1215
    0.8700    0.1153
    0.8800    0.1069
    0.8900    0.0964
    0.9000    0.0860
    0.9100    0.0712
    0.9200    0.0565
    0.9300    0.0420
    0.9400    0.0279
    0.9500    0.0137
    0.9600    0.0078
    0.9700    0.0019
    0.9800    0.0003
    0.9900    0.0026
    1.0000    0.0050
/\relax

}\relax
\endpicture
}
\setbox\figurethree=\vbox{\hsize=\xfiglen
\beginpicture
  \setcoordinatesystem units <\xfiglen,\yfiglen>  point at 0.0 0.0
  \setplotarea x from 0 to 1, y from 0 to 0.17
\scriptsize
  \axis bottom shiftedto y=0.0 ticks short numbered from 0 to 1 by 0.2 /
  \axis left ticks short numbered from 0 to 0.17 by 0.05 /
\footnotesize
\put {$\kappa(x)$} [lt] at 0.02 0.17
\put {$\alpha=0.5$} [lt] at 0.4 0.16
\setlinear
\setsolid
\Red{\relax   

}\relax
\setdashes <2pt>
\Black{\relax 

}\relax
\setplotsymbol ({\sixrm .})
\setsolid
\Blue{\relax   
\plot
         0    0.0000
    0.0100    0.0004
    0.0200    0.0007
    0.0300    0.0018
    0.0400    0.0036
    0.0500    0.0054
    0.0600    0.0097
    0.0700    0.0141
    0.0800    0.0200
    0.0900    0.0274
    0.1000    0.0348
    0.1100    0.0446
    0.1200    0.0544
    0.1300    0.0645
    0.1400    0.0750
    0.1500    0.0855
    0.1600    0.0945
    0.1700    0.1035
    0.1800    0.1107
    0.1900    0.1162
    0.2000    0.1216
    0.2100    0.1225
    0.2200    0.1234
    0.2300    0.1221
    0.2400    0.1185
    0.2500    0.1149
    0.2600    0.1079
    0.2700    0.1010
    0.2800    0.0932
    0.2900    0.0846
    0.3000    0.0760
    0.3100    0.0675
    0.3200    0.0590
    0.3300    0.0514
    0.3400    0.0445
    0.3500    0.0376
    0.3600    0.0330
    0.3700    0.0285
    0.3800    0.0253
    0.3900    0.0235
    0.4000    0.0216
    0.4100    0.0225
    0.4200    0.0233
    0.4300    0.0252
    0.4400    0.0281
    0.4500    0.0310
    0.4600    0.0350
    0.4700    0.0390
    0.4800    0.0430
    0.4900    0.0470
    0.5000    0.0511
    0.5100    0.0544
    0.5200    0.0577
    0.5300    0.0606
    0.5400    0.0630
    0.5500    0.0655
    0.5600    0.0676
    0.5700    0.0698
    0.5800    0.0725
    0.5900    0.0757
    0.6000    0.0789
    0.6100    0.0842
    0.6200    0.0895
    0.6300    0.0959
    0.6400    0.1033
    0.6500    0.1107
    0.6600    0.1187
    0.6700    0.1266
    0.6800    0.1336
    0.6900    0.1397
    0.7000    0.1458
    0.7100    0.1487
    0.7200    0.1516
    0.7300    0.1532
    0.7400    0.1536
    0.7500    0.1539
    0.7600    0.1530
    0.7700    0.1520
    0.7800    0.1505
    0.7900    0.1483
    0.8000    0.1461
    0.8100    0.1426
    0.8200    0.1390
    0.8300    0.1353
    0.8400    0.1312
    0.8500    0.1272
    0.8600    0.1224
    0.8700    0.1175
    0.8800    0.1101
    0.8900    0.1002
    0.9000    0.0904
    0.9100    0.0735
    0.9200    0.0566
    0.9300    0.0403
    0.9400    0.0248
    0.9500    0.0092
    0.9600    0.0056
    0.9700    0.0019
    0.9800    0.0013
    0.9900    0.0038
    1.0000    0.0063
/\relax

}\relax
\endpicture
}
\setbox\figurefour=\vbox{\hsize=\xfiglen
\beginpicture
  \setcoordinatesystem units <\xfiglen,\yfiglen>  point at 0.0 0.0
  \setplotarea x from 0 to 1, y from 0 to 0.17
\scriptsize
  \axis bottom shiftedto y=0.0 ticks short numbered from 0 to 1 by 0.2 /
  \axis left ticks short numbered from 0 to 0.17 by 0.05 /
\footnotesize
\put {$\kappa(x)$} [lt] at 0.02 0.17
\put {$\alpha=0.25$} [lt] at 0.4 0.16
\setlinear
\setsolid
\Red{\relax   

}\relax
\setdashes <2pt>
\Black{\relax 

}\relax
\setplotsymbol ({\sixrm .})
\setsolid
\Blue{\relax   
\plot
         0    0.0000
    0.0100    0.0003
    0.0200    0.0006
    0.0300    0.0015
    0.0400    0.0031
    0.0500    0.0048
    0.0600    0.0091
    0.0700    0.0135
    0.0800    0.0195
    0.0900    0.0272
    0.1000    0.0349
    0.1100    0.0453
    0.1200    0.0556
    0.1300    0.0664
    0.1400    0.0774
    0.1500    0.0885
    0.1600    0.0976
    0.1700    0.1068
    0.1800    0.1140
    0.1900    0.1192
    0.2000    0.1244
    0.2100    0.1246
    0.2200    0.1248
    0.2300    0.1227
    0.2400    0.1183
    0.2500    0.1139
    0.2600    0.1063
    0.2700    0.0988
    0.2800    0.0906
    0.2900    0.0820
    0.3000    0.0733
    0.3100    0.0654
    0.3200    0.0576
    0.3300    0.0506
    0.3400    0.0446
    0.3500    0.0386
    0.3600    0.0346
    0.3700    0.0306
    0.3800    0.0277
    0.3900    0.0258
    0.4000    0.0239
    0.4100    0.0244
    0.4200    0.0249
    0.4300    0.0264
    0.4400    0.0289
    0.4500    0.0314
    0.4600    0.0349
    0.4700    0.0384
    0.4800    0.0421
    0.4900    0.0458
    0.5000    0.0495
    0.5100    0.0529
    0.5200    0.0563
    0.5300    0.0594
    0.5400    0.0621
    0.5500    0.0649
    0.5600    0.0673
    0.5700    0.0697
    0.5800    0.0726
    0.5900    0.0758
    0.6000    0.0790
    0.6100    0.0844
    0.6200    0.0898
    0.6300    0.0962
    0.6400    0.1037
    0.6500    0.1112
    0.6600    0.1191
    0.6700    0.1271
    0.6800    0.1341
    0.6900    0.1402
    0.7000    0.1462
    0.7100    0.1487
    0.7200    0.1512
    0.7300    0.1526
    0.7400    0.1528
    0.7500    0.1530
    0.7600    0.1525
    0.7700    0.1519
    0.7800    0.1507
    0.7900    0.1489
    0.8000    0.1470
    0.8100    0.1431
    0.8200    0.1392
    0.8300    0.1351
    0.8400    0.1307
    0.8500    0.1264
    0.8600    0.1219
    0.8700    0.1173
    0.8800    0.1104
    0.8900    0.1009
    0.9000    0.0915
    0.9100    0.0741
    0.9200    0.0568
    0.9300    0.0401
    0.9400    0.0242
    0.9500    0.0083
    0.9600    0.0051
    0.9700    0.0020
    0.9800    0.0015
    0.9900    0.0041
    1.0000    0.0066
/\relax

}\relax
\endpicture
}
\setbox\figuretwelve=\vbox{\hsize=\xfiglen
\beginpicture
  \setcoordinatesystem units <\xfiglen,\yfiglen>  point at 0.0 0.0
  \setplotarea x from 0 to 1, y from 0 to 0.17
\scriptsize
  \axis bottom shiftedto y=0.0 ticks short numbered from 0 to 1 by 0.2 /
  \axis left ticks short numbered from 0 to 0.17 by 0.05 /
\footnotesize
\put {$\kappa(x)$} [lt] at 0.02 0.17
\put {$\alpha=0.25$} [lt] at 0.4 0.16
\setlinear
\setsolid
\Red{\relax   

}\relax
\setdashes <2pt>
\Black{\relax 

}\relax
\setplotsymbol ({\sixrm .})
\setsolid
\Blue{\relax   
\plot
         0    0.0000
    0.0100    0.0001
    0.0200    0.0002
    0.0300    0.0005
    0.0400    0.0010
    0.0500    0.0015
    0.0600    0.0030
    0.0700    0.0046
    0.0800    0.0070
    0.0900    0.0103
    0.1000    0.0136
    0.1100    0.0191
    0.1200    0.0246
    0.1300    0.0312
    0.1400    0.0388
    0.1500    0.0465
    0.1600    0.0554
    0.1700    0.0644
    0.1800    0.0732
    0.1900    0.0818
    0.2000    0.0905
    0.2100    0.0969
    0.2200    0.1034
    0.2300    0.1079
    0.2400    0.1104
    0.2500    0.1129
    0.2600    0.1108
    0.2700    0.1086
    0.2800    0.1045
    0.2900    0.0984
    0.3000    0.0923
    0.3100    0.0842
    0.3200    0.0760
    0.3300    0.0680
    0.3400    0.0600
    0.3500    0.0519
    0.3600    0.0457
    0.3700    0.0395
    0.3800    0.0347
    0.3900    0.0313
    0.4000    0.0279
    0.4100    0.0274
    0.4200    0.0270
    0.4300    0.0276
    0.4400    0.0294
    0.4500    0.0312
    0.4600    0.0340
    0.4700    0.0368
    0.4800    0.0396
    0.4900    0.0425
    0.5000    0.0454
    0.5100    0.0482
    0.5200    0.0510
    0.5300    0.0537
    0.5400    0.0564
    0.5500    0.0591
    0.5600    0.0621
    0.5700    0.0650
    0.5800    0.0686
    0.5900    0.0728
    0.6000    0.0770
    0.6100    0.0833
    0.6200    0.0895
    0.6300    0.0966
    0.6400    0.1044
    0.6500    0.1123
    0.6600    0.1202
    0.6700    0.1280
    0.6800    0.1348
    0.6900    0.1406
    0.7000    0.1463
    0.7100    0.1487
    0.7200    0.1510
    0.7300    0.1523
    0.7400    0.1524
    0.7500    0.1525
    0.7600    0.1518
    0.7700    0.1510
    0.7800    0.1498
    0.7900    0.1481
    0.8000    0.1463
    0.8100    0.1430
    0.8200    0.1398
    0.8300    0.1361
    0.8400    0.1321
    0.8500    0.1281
    0.8600    0.1229
    0.8700    0.1176
    0.8800    0.1099
    0.8900    0.0995
    0.9000    0.0892
    0.9100    0.0723
    0.9200    0.0554
    0.9300    0.0395
    0.9400    0.0246
    0.9500    0.0097
    0.9600    0.0076
    0.9700    0.0054
    0.9800    0.0050
    0.9900    0.0065
    1.0000    0.0079
/\relax

}\relax
\Green{\relax   
\plot
         0    0.0000
    0.0100    0.0001
    0.0200    0.0001
    0.0300    0.0003
    0.0400    0.0006
    0.0500    0.0009
    0.0600    0.0017
    0.0700    0.0025
    0.0800    0.0037
    0.0900    0.0053
    0.1000    0.0069
    0.1100    0.0094
    0.1200    0.0119
    0.1300    0.0150
    0.1400    0.0186
    0.1500    0.0222
    0.1600    0.0268
    0.1700    0.0315
    0.1800    0.0364
    0.1900    0.0417
    0.2000    0.0470
    0.2100    0.0523
    0.2200    0.0576
    0.2300    0.0626
    0.2400    0.0671
    0.2500    0.0717
    0.2600    0.0747
    0.2700    0.0776
    0.2800    0.0795
    0.2900    0.0801
    0.3000    0.0808
    0.3100    0.0790
    0.3200    0.0772
    0.3300    0.0743
    0.3400    0.0702
    0.3500    0.0662
    0.3600    0.0608
    0.3700    0.0555
    0.3800    0.0501
    0.3900    0.0446
    0.4000    0.0392
    0.4100    0.0350
    0.4200    0.0308
    0.4300    0.0277
    0.4400    0.0258
    0.4500    0.0238
    0.4600    0.0245
    0.4700    0.0251
    0.4800    0.0269
    0.4900    0.0299
    0.5000    0.0328
    0.5100    0.0374
    0.5200    0.0419
    0.5300    0.0468
    0.5400    0.0521
    0.5500    0.0573
    0.5600    0.0627
    0.5700    0.0682
    0.5800    0.0736
    0.5900    0.0790
    0.6000    0.0845
    0.6100    0.0901
    0.6200    0.0958
    0.6300    0.1015
    0.6400    0.1075
    0.6500    0.1134
    0.6600    0.1193
    0.6700    0.1253
    0.6800    0.1309
    0.6900    0.1361
    0.7000    0.1414
    0.7100    0.1450
    0.7200    0.1486
    0.7300    0.1512
    0.7400    0.1527
    0.7500    0.1542
    0.7600    0.1536
    0.7700    0.1530
    0.7800    0.1518
    0.7900    0.1498
    0.8000    0.1478
    0.8100    0.1449
    0.8200    0.1420
    0.8300    0.1384
    0.8400    0.1341
    0.8500    0.1298
    0.8600    0.1226
    0.8700    0.1154
    0.8800    0.1060
    0.8900    0.0942
    0.9000    0.0824
    0.9100    0.0671
    0.9200    0.0518
    0.9300    0.0378
    0.9400    0.0252
    0.9500    0.0126
    0.9600    0.0129
    0.9700    0.0133
    0.9800    0.0142
    0.9900    0.0155
    1.0000    0.0169
/\relax

}\relax
\endpicture
}
\setbox\figurethirteen=\vbox{\hsize=\xfiglen
\beginpicture
  \setcoordinatesystem units <\xfiglen,\yfiglen>  point at 0.0 0.0
  \setplotarea x from 0 to 1, y from 0 to 0.17
\scriptsize
  \axis bottom shiftedto y=0.0 ticks short numbered from 0 to 1 by 0.2 /
  \axis left ticks short numbered from 0 to 0.17 by 0.05 /
\footnotesize
\put {$\kappa(x)$} [lt] at 0.02 0.17
\put {$\alpha=0.9$} [lt] at 0.4 0.16
\setlinear
\setsolid
\Red{\relax   

}\relax
\setdashes <2pt>
\Black{\relax 

}\relax
\setplotsymbol ({\sixrm .})
\setsolid
\Blue{\relax   
\plot
         0    0.0000
    0.0100    0.0001
    0.0200    0.0003
    0.0300    0.0007
    0.0400    0.0013
    0.0500    0.0020
    0.0600    0.0037
    0.0700    0.0053
    0.0800    0.0076
    0.0900    0.0106
    0.1000    0.0136
    0.1100    0.0180
    0.1200    0.0225
    0.1300    0.0276
    0.1400    0.0333
    0.1500    0.0391
    0.1600    0.0457
    0.1700    0.0523
    0.1800    0.0589
    0.1900    0.0655
    0.2000    0.0722
    0.2100    0.0779
    0.2200    0.0836
    0.2300    0.0884
    0.2400    0.0922
    0.2500    0.0961
    0.2600    0.0972
    0.2700    0.0984
    0.2800    0.0981
    0.2900    0.0964
    0.3000    0.0946
    0.3100    0.0902
    0.3200    0.0857
    0.3300    0.0804
    0.3400    0.0742
    0.3500    0.0680
    0.3600    0.0612
    0.3700    0.0545
    0.3800    0.0482
    0.3900    0.0423
    0.4000    0.0364
    0.4100    0.0326
    0.4200    0.0288
    0.4300    0.0262
    0.4400    0.0251
    0.4500    0.0239
    0.4600    0.0253
    0.4700    0.0267
    0.4800    0.0292
    0.4900    0.0326
    0.5000    0.0360
    0.5100    0.0406
    0.5200    0.0451
    0.5300    0.0499
    0.5400    0.0548
    0.5500    0.0597
    0.5600    0.0647
    0.5700    0.0696
    0.5800    0.0746
    0.5900    0.0796
    0.6000    0.0847
    0.6100    0.0901
    0.6200    0.0955
    0.6300    0.1011
    0.6400    0.1070
    0.6500    0.1128
    0.6600    0.1188
    0.6700    0.1248
    0.6800    0.1304
    0.6900    0.1357
    0.7000    0.1409
    0.7100    0.1446
    0.7200    0.1483
    0.7300    0.1509
    0.7400    0.1524
    0.7500    0.1539
    0.7600    0.1534
    0.7700    0.1529
    0.7800    0.1517
    0.7900    0.1497
    0.8000    0.1478
    0.8100    0.1448
    0.8200    0.1419
    0.8300    0.1383
    0.8400    0.1339
    0.8500    0.1296
    0.8600    0.1226
    0.8700    0.1155
    0.8800    0.1063
    0.8900    0.0951
    0.9000    0.0838
    0.9100    0.0692
    0.9200    0.0546
    0.9300    0.0407
    0.9400    0.0277
    0.9500    0.0147
    0.9600    0.0105
    0.9700    0.0063
    0.9800    0.0045
    0.9900    0.0051
    1.0000    0.0056
/\relax

}\relax
\Green{\relax   
\plot
         0    0.0000
    0.0100    0.0001
    0.0200    0.0001
    0.0300    0.0003
    0.0400    0.0006
    0.0500    0.0009
    0.0600    0.0017
    0.0700    0.0025
    0.0800    0.0037
    0.0900    0.0053
    0.1000    0.0069
    0.1100    0.0094
    0.1200    0.0119
    0.1300    0.0150
    0.1400    0.0186
    0.1500    0.0222
    0.1600    0.0268
    0.1700    0.0315
    0.1800    0.0364
    0.1900    0.0417
    0.2000    0.0470
    0.2100    0.0523
    0.2200    0.0576
    0.2300    0.0626
    0.2400    0.0671
    0.2500    0.0717
    0.2600    0.0747
    0.2700    0.0776
    0.2800    0.0795
    0.2900    0.0801
    0.3000    0.0808
    0.3100    0.0790
    0.3200    0.0772
    0.3300    0.0743
    0.3400    0.0702
    0.3500    0.0662
    0.3600    0.0608
    0.3700    0.0555
    0.3800    0.0501
    0.3900    0.0446
    0.4000    0.0392
    0.4100    0.0350
    0.4200    0.0308
    0.4300    0.0277
    0.4400    0.0258
    0.4500    0.0238
    0.4600    0.0245
    0.4700    0.0251
    0.4800    0.0269
    0.4900    0.0299
    0.5000    0.0328
    0.5100    0.0374
    0.5200    0.0419
    0.5300    0.0468
    0.5400    0.0521
    0.5500    0.0573
    0.5600    0.0627
    0.5700    0.0682
    0.5800    0.0736
    0.5900    0.0790
    0.6000    0.0845
    0.6100    0.0901
    0.6200    0.0958
    0.6300    0.1015
    0.6400    0.1075
    0.6500    0.1134
    0.6600    0.1193
    0.6700    0.1253
    0.6800    0.1309
    0.6900    0.1361
    0.7000    0.1414
    0.7100    0.1450
    0.7200    0.1486
    0.7300    0.1512
    0.7400    0.1527
    0.7500    0.1542
    0.7600    0.1536
    0.7700    0.1530
    0.7800    0.1518
    0.7900    0.1498
    0.8000    0.1478
    0.8100    0.1449
    0.8200    0.1420
    0.8300    0.1384
    0.8400    0.1341
    0.8500    0.1298
    0.8600    0.1226
    0.8700    0.1154
    0.8800    0.1060
    0.8900    0.0942
    0.9000    0.0824
    0.9100    0.0671
    0.9200    0.0518
    0.9300    0.0378
    0.9400    0.0252
    0.9500    0.0126
    0.9600    0.0129
    0.9700    0.0133
    0.9800    0.0142
    0.9900    0.0155
    1.0000    0.0169
/\relax

}\relax
\endpicture
}
\xfiglen=4true in
\yfiglen=2true in
\setbox\figurefive=\vbox{\hsize=\xfiglen
\beginpicture
  \setcoordinatesystem units <0.02\xfiglen,0.06\yfiglen>  point at 0.0 0.0
  \setplotarea x from 0 to 30, y from 0 to -12
\scriptsize
  \axis bottom shiftedto y=-12 ticks short numbered from 0 to 30 by 5 /
  \axis left ticks short numbered from -12 to 0 by 4 /
\footnotesize
\put {\copy\figurelegendone} [lb] at 3 -10.5
\scriptsize
\put {$c=1$} [lb] at 13 -10.5
\put {$\log_{10}(\sigma_n)$} [l] at 2 0.0
\put {$n$} [rb] at 30 -11.6
\multiput {$\star$} at
    1.0000    0.0373
    2.0000   -0.8989
    3.0000   -1.4315
    4.0000   -1.7568
    5.0000   -1.9959
    6.0000   -2.1860
    7.0000   -2.3446
    8.0000   -2.4814
    9.0000   -2.6021
   10.0000   -2.7106
   11.0000   -2.8097
   12.0000   -2.9011
   13.0000   -2.9863
   14.0000   -3.0663
   15.0000   -3.1420
   16.0000   -3.2141
   17.0000   -3.2830
   18.0000   -3.3492
   19.0000   -3.4131
   20.0000   -3.4748
   21.0000   -3.5345
   22.0000   -3.5922
   23.0000   -3.6479
   24.0000   -3.7012
   25.0000   -3.7514
   26.0000   -3.7975
   27.0000   -3.8413
   28.0000   -3.8888
   29.0000   -3.9432
   30.0000   -4.0054
/
\multiput {$\diamond$} at
    1.0000    0.0367
    2.0000   -0.9021
    3.0000   -1.4400
    4.0000   -1.7725
    5.0000   -2.0197
    6.0000   -2.2189
    7.0000   -2.3874
    8.0000   -2.5347
    9.0000   -2.6666
   10.0000   -2.7870
   11.0000   -2.8984
   12.0000   -3.0026
   13.0000   -3.1011
   14.0000   -3.1949
   15.0000   -3.2848
   16.0000   -3.3713
   17.0000   -3.4550
   18.0000   -3.5361
   19.0000   -3.6148
   20.0000   -3.6909
   21.0000   -3.7641
   22.0000   -3.8347
   23.0000   -3.9047
   24.0000   -3.9772
   25.0000   -4.0541
   26.0000   -4.1362
   27.0000   -4.2237
   28.0000   -4.3172
   29.0000   -4.4174
   30.0000   -4.5254
/
\multiput {$\circ$} at
    1.0000    0.0358
    2.0000   -0.9088
    3.0000   -1.4614
    4.0000   -1.8200
    5.0000   -2.1024
    6.0000   -2.3446
    7.0000   -2.5631
    8.0000   -2.7668
    9.0000   -2.9609
   10.0000   -3.1488
   11.0000   -3.3326
   12.0000   -3.5135
   13.0000   -3.6921
   14.0000   -3.8692
   15.0000   -4.0469
   16.0000   -4.2295
   17.0000   -4.4213
   18.0000   -4.6251
   19.0000   -4.8419
   20.0000   -5.0725
   21.0000   -5.3175
   22.0000   -5.5778
   23.0000   -5.8546
   24.0000   -6.1493
   25.0000   -6.4631
   26.0000   -6.7978
   27.0000   -7.1557
   28.0000   -7.5387
   29.0000   -7.9495
   30.0000   -8.3924
/
\multiput {$\bullet$} at
    1.0000    0.0355
    2.0000   -0.9105
    3.0000   -1.4677
    4.0000   -1.8364
    5.0000   -2.1336
    6.0000   -2.3952
    7.0000   -2.6373
    8.0000   -2.8687
    9.0000   -3.0943
   10.0000   -3.3172
   11.0000   -3.5392
   12.0000   -3.7611
   13.0000   -3.9844
   14.0000   -4.2122
   15.0000   -4.4496
   16.0000   -4.7013
   17.0000   -4.9701
   18.0000   -5.2569
   19.0000   -5.5627
   20.0000   -5.8886
   21.0000   -6.2356
   22.0000   -6.6049
   23.0000   -6.9980
   24.0000   -7.4168
   25.0000   -7.8627
   26.0000   -8.3372
   27.0000   -8.8449
   28.0000   -9.3867
   29.0000   -9.9623
   30.0000  -10.5896
/
\endpicture
}
\xfiglen=3true in
\yfiglen=2true in
\setbox\figuresix=\vbox{\hsize=\xfiglen
\beginpicture
  \setcoordinatesystem units <0.13\xfiglen,0.018\yfiglen>  point at 0.0 0.0
  \setplotarea x from -6.2 to 0, y from 0 to 40
\scriptsize
  \axis bottom shiftedto y=0 ticks short numbered from -6 to 0 by 1 /
  \axis right shiftedto y=0 ticks short numbered from 0 to 40 by 5 /
\footnotesize
\multiput {$\star$} at  
   -0.1967    3.3386
   -0.5557    6.8409
   -1.0200   10.4501
   -1.5689   14.1461
   -2.1905   17.9171
   -2.8768   21.7552
   -3.6218   25.6548
   -4.4209   29.6113
   -5.2703   33.6211
   -6.1670   37.6812
/
\multiput {$\circ$} at  
   -0.1110    3.2527
   -0.2221    6.5054
   -0.3331    9.7581
   -0.4441   13.0108
   -0.5552   16.2635
   -0.6662   19.5162
   -0.7773   22.7689
   -0.8883   26.0216
   -0.9993   29.2743
   -1.1104   32.5270
   -1.2214   35.7797
   -1.3324   39.0324
/
\multiput {$\diamond$} at  
   -0.0627    3.2043
   -0.0886    6.3718
   -0.1085    9.5333
   -0.1253   12.6917
   -0.1401   15.8481
   -0.1535   19.0031
   -0.1658   22.1569
   -0.1772   25.3100
   -0.1880   28.4623
   -0.1982   31.6141
   -0.2078   34.7654
   -0.2171   37.9162
/
\scriptsize
\put {$\alpha=\frac{1}{2}$, $\;c=1$} [l] at -6 15
\put {\copy\figurelegendtwo} [lb] at -6 2
\endpicture
}
\setbox\figureseven=\vbox{\hsize=\xfiglen
\beginpicture
  \setcoordinatesystem units <0.029\xfiglen,0.008\yfiglen>  point at 0.0 0.0
  \setplotarea x from -25 to 0, y from 0 to 90
\scriptsize
  \axis bottom shiftedto y=0 ticks short numbered from -25 to 0 by 5 /
  \axis right shiftedto y=0 ticks short numbered from 0 to 90 by 10 /
\footnotesize
\multiput {$\star$} at  
   -0.4733    0.8073
   -1.2033    2.0345
   -2.0814    3.4938
   -3.0744    5.1279
   -4.1642    6.9057
   -5.3393    8.8074
   -6.5914   10.8187
   -7.9143   12.9290
   -9.3030   15.1301
  -10.7535   17.4150
/
\multiput {$\circ$} at  
   -1.3612    2.5378
   -3.2853    6.3631
   -5.4497   10.8680
   -7.7657   15.8630
  -10.1897   21.2449
  -12.6964   26.9473
  -15.2695   32.9233
  -17.8982   39.1380
  -20.5745   45.5644
  -23.2927   52.1809
/
\multiput {$\diamond$} at  
   -1.8611    6.5699
   -3.4835   14.7744
   -5.0047   23.4568
   -6.4928   32.4216
   -7.9765   41.5885
   -9.4705   50.9140
  -10.9829   60.3713
  -12.5184   69.9422
  -14.0801   79.6134
  -15.6696   89.3749
/
\scriptsize
\put {$\alpha=\frac{1}{2}$, $\;c=5$} [l] at -24.9 30
\put {\copy\figurelegendtwo} [lb] at -24.8 2
\endpicture
}
\setbox\figureeight=\vbox{\hsize=\xfiglen
\beginpicture
  \setcoordinatesystem units <0.0278\xfiglen,0.024\yfiglen>  point at 0.0 0.0
  \setplotarea x from -30 to 0, y from 0 to 32
\scriptsize
  \axis bottom shiftedto y=0 ticks short numbered from -25 to 0 by 5 /
  \axis right shiftedto y=0 ticks short numbered from 0 to 32 by 5 /
\footnotesize
\multiput {$\star$} at  
   -0.4430    3.1872
   -1.6858    6.3792
   -3.7223    9.4797
   -6.5929   12.4052
  -10.3760   15.0823
  -15.1979   17.4609
  -21.2370   19.5516
  -28.6841   21.4866
/
\multiput {$\circ$} at  
   -0.2478    3.1730
   -0.6562    6.3599
   -1.1607    9.5522
   -1.7401   12.7473
   -2.3828   15.9439
   -3.0809   19.1412
   -3.8291   22.3385
   -4.6229   25.5354
   -5.4591   28.7316
   -6.3350   31.9268
/
\multiput {$\diamond$} at  
   -0.1392    3.1611
   -0.2596    6.3200
   -0.3738    9.4780
   -0.4842   12.6355
   -0.5918   15.7927
   -0.6973   18.9495
   -0.8010   22.1062
   -0.9033   25.2626
   -1.0042   28.4189
   -1.1040   31.5751
/
\scriptsize
\put {$\alpha=\frac{9}{10}$, $\;c=1$} [lb] at -29.7 11
\put {\copy\figurelegendtwo} [lb] at -29 2
\endpicture
}
\setbox\figurenine=\vbox{\hsize=\xfiglen
\beginpicture
  \setcoordinatesystem units <0.025\xfiglen,0.0047\yfiglen>  point at 0.0 0.0
  \setplotarea x from -30 to 0, y from 0 to 160
\scriptsize
  \axis bottom shiftedto y=0 ticks short numbered from -30 to 0 by 5 /
  \axis right shiftedto y=0 ticks short numbered from 0 to 160 by 40 /
\footnotesize
\multiput {$\star$} at  
   -0.3712   15.7632
   -1.3890   31.6111
   -3.0090   47.5233
   -5.2110   63.4845
   -7.9831   79.4819
  -11.3170   95.5034
  -15.2074  111.5380
  -19.6505  127.5753
  -24.6437  143.6052
  -30.1857  159.6179
/
\multiput {$\circ$} at  
   -0.2091   15.7400
   -0.5522   31.4994
   -0.9746   47.2700
   -1.4584   63.0491
   -1.9938   78.8350
   -2.5742   94.6269
   -3.1949  110.4241
   -3.8525  126.2259
   -4.5440  142.0321
   -5.2671  157.8422
/
\multiput {$\diamond$} at  
   -0.1179   15.7263
   -0.2200   31.4501
   -0.3169   47.1732
   -0.4105   62.8958
   -0.5018   78.6180
   -0.5913   94.3399
   -0.6793  110.0616
   -0.7660  125.7831
   -0.8517  141.5045
   -0.9364  157.2257
/
\scriptsize
\put {$\alpha=\frac{9}{10}$, $\;c=5$} [lb] at -30 50
\put {\copy\figurelegendtwo} [lb] at -29 2
\endpicture
}
\setbox\figureten=\vbox{\hsize=\xfiglen
\beginpicture
  \setcoordinatesystem units <0.009\xfiglen,0.013\yfiglen>  point at 0.0 0.0
  \setplotarea x from -80 to 0, y from 0 to 60
\scriptsize
  \axis bottom shiftedto y=0 ticks short numbered from -80 to 0 by 20 /
  \axis right shiftedto y=0 ticks short numbered from 0 to 60 by 10 /
\footnotesize
\multiput {$\star$} at  
   -0.3937    6.5322
   -1.2549   13.1640
   -2.4725   19.7357
   -3.9638   26.1725
   -5.6554   32.4400
   -7.4864   38.5286
   -9.4097   44.4433
  -11.3903   50.1956
  -13.4028   55.7999
  -15.4290   61.2705
/
\multiput {$\circ$} at  
   -0.9859    6.8911
   -3.1132   13.9700
   -6.0368   20.8379
   -9.5052   27.3574
  -13.3300   33.5031
  -17.3825   39.2972
  -21.5785   44.7776
  -25.8633   49.9833
  -30.2016   54.9498
  -34.5703   59.7072
/
\multiput {$\diamond$} at  
   -1.9717    7.4641
   -6.1345   15.1739
  -11.6541   22.3859
  -17.9896   28.9578
  -24.8107   34.9359
  -31.9276   40.4090
  -39.2307   45.4617
  -46.6555   50.1640
  -54.1619   54.5714
  -61.7249   58.7279
/
\multiput {$\bullet$} at  
   -3.9015    8.5621
  -11.7857   17.3634
  -21.7017   25.1794
  -32.6488   32.0097
  -44.1485   38.0511
  -55.9675   43.4809
  -67.9843   48.4313
  -80.1310   52.9980
/
\put {\copy\figurelegendthree} [lb] at -95 5
\endpicture
}
\setbox\figureeleven=\vbox{\hsize=\xfiglen
\beginpicture
  \setcoordinatesystem units <0.095\xfiglen,0.0066\yfiglen>  point at 0.0 0.0
  \setplotarea x from -7.5 to 0, y from 0 to 120
\scriptsize
  \axis bottom shiftedto y=0 ticks short numbered from -7 to 0 by 1 /
  \axis right shiftedto y=0 ticks short numbered from 0 to 120 by 20 /
\footnotesize
\multiput {$\star$} at  
   -0.3941    3.3054
   -1.2038    7.1645
   -1.8993   11.5278
   -2.3449   16.0617
   -2.6320   20.6148
   -2.8313   25.1546
   -2.9791   29.6780
   -3.0945   34.1873
   -3.1883   38.6854
   -3.2666   43.1747
   -3.3337   47.6569
   -3.3922   52.1337
   -3.4440   56.6059
   -3.4904   61.0744
   -3.5324   65.5399
   -3.5707   70.0028
   -3.6060   74.4635
   -3.6387   78.9225
   -3.6691   83.3798
   -3.6975   87.8358
   -3.7242   92.2906
   -3.7494   96.7443
   -3.7732  101.1972
   -3.7958  105.6492
   -3.8173  110.1005
   -3.8378  114.5512
   -3.8574  119.0013
/
\multiput {$\circ$} at  
   -2.5123    4.3453
   -4.4717   12.7661
   -5.0481   20.8981
   -5.3647   28.8273
   -5.5813   36.6677
   -5.7458   44.4616
   -5.8786   52.2280
   -5.9901   59.9769
   -6.0862   67.7138
   -6.1708   75.4423
   -6.2465   83.1645
   -6.3149   90.8820
   -6.3774   98.5958
   -6.4349  106.3068
   -6.4883  114.0154
/
\multiput {$\diamond$} at  
   -4.3113    7.1883
   -5.5001   18.5196
   -5.9480   29.2715
   -6.2248   39.8619
   -6.4261   50.3849
   -6.5847   60.8732
   -6.7160   71.3413
   -6.8282   81.7967
   -6.9263   92.2434
   -7.0135  102.6841
   -7.0921  113.1204
/
\multiput {$\bullet$} at  
   -5.3252   11.6667
   -6.1887   26.7393
   -6.5847   41.3710
   -6.8457   55.8854
   -7.0416   70.3515
   -7.1992   84.7930
   -7.3313   99.2203
   -7.4453  113.6387
/
\endpicture
}
\xfiglen=4true in
\yfiglen=2true in
\setbox\figurefourteen=\vbox{\hsize=\xfiglen
\beginpicture
  \setcoordinatesystem units <0.02\xfiglen,0.06\yfiglen>  point at 0.0 0.0
  \setplotarea x from 0 to 30, y from 0 to -12
\scriptsize
  \axis bottom shiftedto y=-12 ticks short numbered from 0 to 30 by 5 /
  \axis left ticks short numbered from -12 to 0 by 4 /
\footnotesize
\put {\copy\figurelegendone} [lb] at 2 -10.5
\scriptsize
\put {$c=5$} [lb] at 13 -10.5
\scriptsize
\put {$\log_{10}(\sigma_n)$} [l] at 1 0.0
\put {$n$} [rb] at 30 -11.6
\multiput {$\bullet$} at  
    1.0000   -1.0330
    2.0000   -2.2119
    3.0000   -2.8233
    4.0000   -3.2073
    5.0000   -3.5151
    6.0000   -3.7693
    7.0000   -3.9902
    8.0000   -4.1805
    9.0000   -4.3484
   10.0000   -4.4984
   11.0000   -4.6368
   12.0000   -4.7668
   13.0000   -4.8912
   14.0000   -5.0115
   15.0000   -5.1287
   16.0000   -5.2433
   17.0000   -5.3556
   18.0000   -5.4661
   19.0000   -5.5755
   20.0000   -5.6867
   21.0000   -5.8024
   22.0000   -5.9250
   23.0000   -6.0551
   24.0000   -6.1934
   25.0000   -6.3406
   26.0000   -6.4973
   27.0000   -6.6647
   28.0000   -6.8441
   29.0000   -7.0371
   30.0000   -7.2457
/
\multiput {$\circ$} at  
    1.0000   -1.0329
    2.0000   -2.2088
    3.0000   -2.8107
    4.0000   -3.1801
    5.0000   -3.4718
    6.0000   -3.7107
    7.0000   -3.9197
    8.0000   -4.1020
    9.0000   -4.2651
   10.0000   -4.4105
   11.0000   -4.5427
   12.0000   -4.6638
   13.0000   -4.7768
   14.0000   -4.8835
   15.0000   -4.9855
   16.0000   -5.0838
   17.0000   -5.1791
   18.0000   -5.2719
   19.0000   -5.3621
   20.0000   -5.4502
   21.0000   -5.5362
   22.0000   -5.6223
   23.0000   -5.7108
   24.0000   -5.8041
   25.0000   -5.9030
   26.0000   -6.0083
   27.0000   -6.1204
   28.0000   -6.2402
   29.0000   -6.3686
   30.0000   -6.5070
/
\multiput {$\diamond$} at  
    1.0000   -1.1831
    2.0000   -2.3521
    3.0000   -2.9333
    4.0000   -3.2705
    5.0000   -3.5231
    6.0000   -3.7202
    7.0000   -3.8877
    8.0000   -4.0306
    9.0000   -4.1588
   10.0000   -4.2732
   11.0000   -4.3790
   12.0000   -4.4761
   13.0000   -4.5674
   14.0000   -4.6529
   15.0000   -4.7344
   16.0000   -4.8116
   17.0000   -4.8858
   18.0000   -4.9568
   19.0000   -5.0255
   20.0000   -5.0917
   21.0000   -5.1559
   22.0000   -5.2182
   23.0000   -5.2787
   24.0000   -5.3374
   25.0000   -5.3940
   26.0000   -5.4487
   27.0000   -5.5005
   28.0000   -5.5489
   29.0000   -5.5936
   30.0000   -5.6410
/
\multiput {$\star$} at  
    1.0000   -1.0326
    2.0000   -2.2002
    3.0000   -2.7780
    4.0000   -3.1105
    5.0000   -3.3577
    6.0000   -3.5489
    7.0000   -3.7099
    8.0000   -3.8461
    9.0000   -3.9673
   10.0000   -4.0746
   11.0000   -4.1730
   12.0000   -4.2628
   13.0000   -4.3468
   14.0000   -4.4249
   15.0000   -4.4991
   16.0000   -4.5690
   17.0000   -4.6361
   18.0000   -4.7001
   19.0000   -4.7618
   20.0000   -4.8212
   21.0000   -4.8788
   22.0000   -4.9344
   23.0000   -4.9882
   24.0000   -5.0402
   25.0000   -5.0900
   26.0000   -5.1371
   27.0000   -5.1804
   28.0000   -5.2182
   29.0000   -5.2555
   30.0000   -5.3021
/
\endpicture
}
\xfiglen=2.4 true in
\yfiglen=1.25 true in
\setbox\figurefifteen=\vbox{\hsize=\xfiglen
\beginpicture
  \setcoordinatesystem units <0.5\xfiglen,\yfiglen>  point at -0.65 -0.15
  \setplotarea x from -0.65 to 1, y from -0.15 to 1
\scriptsize
  \axis bottom shiftedto y=0.0 / 
  \axis left shiftedto y=0 / 
\footnotesize
\setlinear
\setdots <2pt>
\put {$x$} [rt]  at 1 -0.04
\put {$y$} [lt]  at 0.02 1
\put {$\bullet$}  at 0.5 0.866
\put {$s$} [lb]  at 0.52 0.89
\plot 0 0 0.5 0.866 /
\setsolid
\setplotsymbol ({\sevenrm .})
\OliveGreen{\relax
\arrow <10pt> [0.3, 0.6] from 0 0 to 0.5 0.5 \relax}\relax  
\put {$b\lambda^\beta s^\alpha$} [lb]   at 0.52 0.52
\OliveGreen{\relax
\arrow <10pt> [0.3, 0.6] from 0 0 to -0.5 0.866 \relax}\relax  
\put {$s^2$} [rb]  at -0.46 0.89
\OliveGreen{\relax
\arrow <10pt> [0.3, 0.6] from 0 0 to 0.8 0 \relax}\relax 
\put {$c^2\lambda$} [lb]  at 0.8 0.05

\endpicture
}

\subsection{Location of the poles of the relaxation functions}

Motivated by their role for the degree of ill-posedness of the inverse problem, we develop some further results -- beyond those stated in Lemma \oldref{lem:what} as essentials for our uniqueness proof -- for each of the two models under consideration and also provide some computational results with plots of these poles for several parameter configurations.

{\bf The {\sc cwch} model:}\quad
Poles of the {\sc cwch} model are the roots of the function
$$ \omega_{\text{\thinspace\sc cwch}}(s) :=  s^2 + b\lambda^\beta s^\alpha +  c^2\lambda = 0\,.$$
We recall that according to \cite[Lemma 4.1]{KaltenbacherRundell:2021b} which relies on Rouch\'e's theorem, they lie in the left hand complex plane and this is easily shown by the following alternative argument.

\setlength{\columnsep}{3pt}%
\begin{wrapfigure}[9]{r}{0.4\textwidth}
\centering
    {\copy\figurefifteen}
\end{wrapfigure}
Suppose $s$ is a root in the first quadrant.
Then ${\bf s}$, the line joining the origin to the point $s$ can be
split into a component in the direction of the positive real axis and one
in the direction of the positive imaginary axis.
Then since $\alpha\leq 1$, ${\bf{s}^\alpha}$ has components in the same
directions.
Similarly, the vector ${\bf{s}^{2}}$ has a component parallel to
the real axis and again one in the direction of the positive imaginary axis.
Since $\lambda\geq 0$, the same is true of the vector $b\lambda^\beta {\bf s}^{\alpha}$.
The third vector representing $c^2\lambda {\bf x}$ points along the real axis.
However, the sum of these three vectors cannot add to zero
contradicting the claim the root $s$ lay in the first quadrant.
An identical argument shows $s$ cannot lie in the fourth quadrant and
hence cannot lie in the right half plane.



For $b=0$ the poles are along the imaginary axis and spaced exactly
as the eigenvalue sequence $\{\lambda_n\}$ stretched by the factor $c^2$.
As $b$ increases so does the (negative) real component of the poles which
follow a curve whose rough slope is determined by the ratio of $b$ and $c^2$.
The powers $\alpha$ and $\beta$ also are a factor that 
influence the skewness of the curve along which the poles align.
The magnitude of the real and imaginary parts 
show the relative strengths of the damping and oscillation effects
respectively in the equation.

\begin{figure}[ht]
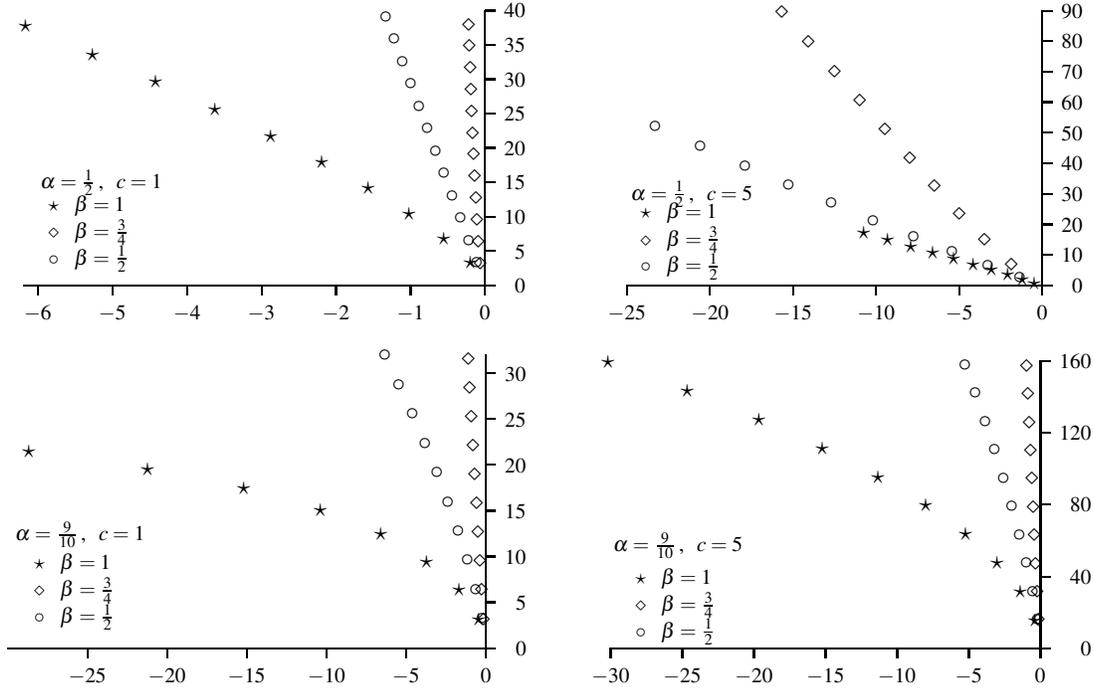

\hbox to \hsize{\hss\copy\figuresix\hss\hss \copy\figureseven\hss}
\vskip10pt
\hbox to \hsize{\hss\copy\figureeight\hss\hss \copy\figurenine\hss}
\caption{
{\small {}\bf Roots of $\;\omega_{\text{\thinspace\sc cwch}}(s)\;$ for various $\alpha$, $\beta$, $c$
values.}}
\label{fig:roots_CWCH}
\end{figure}

The roots of $\omega_{\text{\thinspace\sc cwch}}(s)$
are shown in Figure~\oldref{fig:roots_CWCH} with $b=0.1$, 
 $\lambda_n = n^2\pi^2$, and for both $c=1$ and $c=5$, as well as
$\alpha=\frac12$ and $\alpha=\frac{9}{10}$ and illustrate the above point.

Some notes on how these poles were computed.
For rational $\alpha = p/q$, $\omega_{\text{\thinspace\sc cwch}}(s)$ can be written as
$z^{2q} + B z^p + C$  with $B=b\lambda^\beta$, $C = c^2\lambda$
and where $s = z^p$.
Now the $2q^{\rm th}$ degree polynomial can be represented as the
characteristic polynomial of a $2q\times 2q$ matrix.
Then the roots of this polynomial are calculated by computing the
eigenvalues of the companion matrix.
This gives a good approximation even for reasonably large $q$ values but 
additional care must be taken,
see, for example, \cite{EdelmanMurakami:1995}.
Given now the values of $\{z_n\}$ for $\lambda\in \{\lambda_n\}$, one can recover $\{s_n\}$ from
$s_n = z_n^p$.
This is subject to considerable round-off error for even modest values of $p$.
However it is usually sufficient as an initial approximation for
Newton's method to then compute a more exact value of the roots of $\omega$
to desired accuracy.
This is also successful for real $\alpha$ by first taking a rational
approximation $\alpha \approx p/q$ for the initial approximation
of the roots and then proceeding as above.

{\bf The {\sc fz} model:}\quad
Poles of the fractional Zener model are the roots of the function
$$
\omega_{\,\text{\sc fz}}(s) := b_2 s^{2+\alpha_2} + s^2 + b_1\lambda s^{\alpha_1} +  c^2\lambda = 0.
$$
There is a more complex relationship here and more constants whose value
can affect the outcome.
In the case that $\alpha_1=\alpha_2=\alpha$, we can re-write this as 
$$
\omega_{\,\text{\sc fz}}(s) := (b_2 s^\alpha +1)(s^2 + c^2\lambda) + \delta\lambda s^\alpha  = 0,
\qquad \mbox{where }
\delta := b_1-c^2b_2 
$$
and $\delta$ needs to be nonnegative, cf.
\cite[Section {\sc iii.b}]{HolmNaesholm:2011}.
If $\delta=0$ then $\omega_{\,\text{\sc fz}}(s)$ factors.
There will be two roots at $\pm i\,\sqrt{\lambda}\,c$ on the imaginary axis and a 
potential root coming from $s^\alpha + 1/b_2 =0$. 
The latter only exists in case $\alpha=1$ for otherwise writing $s=re^{i\theta}$ with $\theta\in(-\pi,\pi]$ we have that $\alpha\theta\in(-\pi,\pi)$ and therefore $\Im(s^{\alpha}+1)=r^\alpha\sin(\alpha\theta)=0$ implies $\alpha\theta=0$, hence $\Re(s^{\alpha}+1/b_2)=r^\alpha\cos(\alpha\theta)+1/b_2>0$.
In case $\alpha=1$ we obviously have a root at $-1/b_2$, whose modulus, notably, does not increase with $\lambda$, as opposed to the two other complex conjugate roots of $\omega_{\,\text{\sc fz}}$.

Clearly, physical reasoning leads us to the conclusion that in case of a nonegative diffusivity of sound $\delta\geq0$, all poles need to  have nonpositive real part.
However, the complex analysis arguments from \cite[Lemma 4.1]{KaltenbacherRundell:2021b} via Rouch\'e's Theorem, using as a bounding function the dominant power part $f(z) = b_2z^{2+\alpha_2}+c^2\lambda$, does not seem to directly carry over to the {\sc FZ} case. This is basically due to the fact that we cannot say anything about the number of roots of the non-polynomial function $f$.
Additionally, asymptotics in terms of powers of $s$ will be much less effective
here since, for small $b_2$ and/or $\alpha_2$,
the term $s^2$ will be de facto dominant even for relatively large magnitudes
of $s$.

Therefore we have to take a different path to conclude that also in the {\sc fz} case, the poles lie in the left hand complex plane. We do so by means of energy estimates similar to those in section \oldref{sec:analysisFZ}, which basically corresponds to the mentioned physical argument.
As a (partial) counterpart to \cite[Lemma 4.1]{KaltenbacherRundell:2021b} we state the following.
\begin{lemma}
The roots of $\omega_{\,\text{\sc fz}}(s)$ with $\alpha_1=\alpha_2=\alpha$ and $\delta := b_1-c^2b_2 \geq0$ lie in the left hand complex plane.
\end{lemma}
\begin{proof}
We consider the following initial value problem for the relaxation equation
\begin{equation}\label{eqn:ivpw}
b_2 \partial_t^{2+\alpha_2} w + w'' + b_1 \lambda \partial_t^{\alpha_1} c^2\lambda w =0\,, \quad w(0)=0\,, \ w'(0)=1\,, \quad w''(0)=0\,.
\end{equation}
The Laplace transform $\hat{w}$ of its solution satisfies
$$
(b_2 s^{2+\alpha_2} + s^2 + b_1 \lambda s^{\alpha_1} c^2\lambda) \hat{w}(s) = b_2 s^{\alpha_2} + 1
$$
and therefore $\hat{w}(s)=\frac{b_2 s^{\alpha_2} + 1}{\omega_{\,\text{\sc fz}}(s)}$.
Now if $\alpha_1=\alpha_2=\alpha$ and $\delta := b_1-c^2b_2 \geq0$, analogously to the proof of Proposition \oldref{prop:FZalpha12}, we obtain an energy estimate for $\tilde{w}:=b_2 \partial_t^\alpha w + w$ by multiplying \eqref{eqn:ivpw} with $w''$ and integrating with respect to time
$$
\begin{aligned}
&\frac12 |\tilde{w}_t(t)|^2 + \frac{c^2\lambda}{2}|\tilde{w}(t)|^2
+\delta\lambda\Bigl(\frac{b_2}{2}|\partial_t^\alpha w(t)|^2+\frac{1}{2\Gamma(\alpha)t^{1-\alpha}} \int_0^t |\partial_t^\alpha w(\tau)|^2\, d\tau \Bigr)\\
&\leq \frac12 |\tilde{w}_t(0)|^2 + \frac{c^2\lambda}{2}|\tilde{w}(0)|^2
\end{aligned}
$$
for all $t\geq0$. 
This implies uniform boundedness $|\tilde{w}(t)|\leq\sqrt{(c^2\lambda_{\min})^{-1}|\tilde{w}_t(0)|^2 + |\tilde{w}(0)|^2}=:C$ by a constant independent of $\lambda$.
Taking Laplace transforms 
$$
|\hat{\tilde{w}}(s)|=\left|\int_0^\infty e^{-st}\tilde{w}(t)\, dt\right|
\leq \int_0^\infty e^{-\Re(s)t}\, dt \ C = \frac{C}{\Re(s)} 
$$
for $\Re(s)>0$, we see that $\hat{\tilde{w}}(s)$ cannot have any poles in the right half plane.
Due to the identity $\hat{\tilde{w}}(s)=(b_2 s^\alpha +1) \hat{w}(s) -  b_2 s^{\alpha-1} w(0) = (b_2 s^\alpha +1) \hat{w}(s) = \frac{(b_2 s^\alpha + 1)^2}{\omega_{\,\text{\sc fz}}(s)}$ (where the numerator has no zeros in case $\alpha\in(0,1)$), the assertion follows.
\end{proof}

The effect of $\delta$ on the poles in the {\sc fz} model can also be assessed by means of the implicit function theorem, applied to the function
$$
f(r,\theta;\delta)=\left(\begin{array}{l}
b_2r^{2+\alpha_2}\cos((2+\alpha_2)\theta)+r^2\cos(2\theta)+b_1\lambda r^{\alpha_1}\cos(\alpha_1\theta)
+c^2\lambda\\
b_2r^{2+\alpha_2}\sin((2+\alpha_2)\theta)+r^2\sin(2\theta)+b_1\lambda r^{\alpha_1}\sin(\alpha_1\theta)
\end{array}\right)
$$
whose zeros are the magnitudes and arguments of the roots  $s=re^{i\theta}$ of
$\omega^{\,\text{\sc fz}}$.
Now 
$$
\frac{\partial f}{\partial (r,\theta)}=\left(\begin{array}{cc}A_1&A_2\\B_1&B_2\end{array}\right), \quad
\frac{\partial f}{\partial \delta}=\left(\begin{array}{c}C_1\\C_2\end{array}\right)
$$
and using Cramer's rule this yields 
$$
\frac{\partial r}{\partial \delta}= \frac{B_1C_2-B_2C_1}{A_1B_2-A_2B_1}\,, \quad
\frac{\partial \theta}{\partial \delta}= \frac{C_1A_2-C_2A_1}{A_1B_2-A_2B_1}\,,
$$
where $b_1=b_2c^2+\delta$ and for $(r,\theta)$ satisfying
$f(r,\theta;\delta)=0$.
Then
\vskip-6pt
$$
\begin{aligned}
A_1&= (2+\alpha_2)b_2r^{1+\alpha_2}\cos((2+\alpha_2)\theta)+2r\cos(2\theta)+\alpha_1b_1\lambda r^{\alpha_1-1}\cos(\alpha_1\theta)\\
&= -\frac{1}{r}\Bigl(\alpha_2 r^2\cos(2\theta)+(2+\alpha_2-\alpha_1)b_1\lambda r^{\alpha_1}\cos(\alpha_1\theta)+ (2+\alpha_2)c^2\lambda\Bigr)\\
A_2&= (2+\alpha_2)b_2r^{1+\alpha_2}\sin((2+\alpha_2)\theta)+2r\sin(2\theta)+\alpha_1b_1\lambda r^{\alpha_1-1}\sin(\alpha_1\theta)\\
&= -\frac{1}{r}\Bigl(\alpha_2 r^2\sin(2\theta)+(2+\alpha_2-\alpha_1)b_1\lambda r^{\alpha_1}\sin(\alpha_1\theta)\Bigr)\\
B_1&=-r A,\qquad
B_2=rA_1,\qquad
C_1=\lambda r^{\alpha_1}\cos(\alpha_1\theta),\qquad
C_2=\lambda r^{\alpha_1}\sin(\alpha_1\theta)\,.
\end{aligned}
$$
This results in 
$$
\begin{aligned}
A_1B_2-A_2B_1&= r(A_1^2+A_2^2)\\
B_1C_2-B_2C_1&= \lambda r^{\alpha_1} \bigl(\alpha_2 r^2\cos((2\!-\!\alpha_1)\theta) + (2\!+\!\alpha_2\!-\!\alpha_1)b_1\lambda r^{\alpha_1}+(2\!+\!\alpha_2)c^2\lambda \cos(\alpha_1\theta)\bigr)\\
C_1A_2-C_2A_1&= \lambda r^{\alpha_1-1} \bigl(-\alpha_2 r^2\sin((2-\alpha_1)\theta) + (2+\alpha_2)c^2\lambda \sin(\alpha_1\theta)\bigr)\,.
\end{aligned}
$$
and therefore 
$\, \frac{\partial r}{\partial \delta}>0\,$,
$\quad \frac{\partial \theta}{\partial \delta}>0\,$,
for the case of the known roots $r=c\sqrt{\lambda}$, $\theta=\pm\pi/2$ at $\delta=0$, for $\alpha_1=\alpha_2$.
That is, increasing $\delta$ tends to move the poles into the left hand complex plane, which is intuitive in view of its physical role as a diffusivity of sound. 

\begin{figure}[ht]
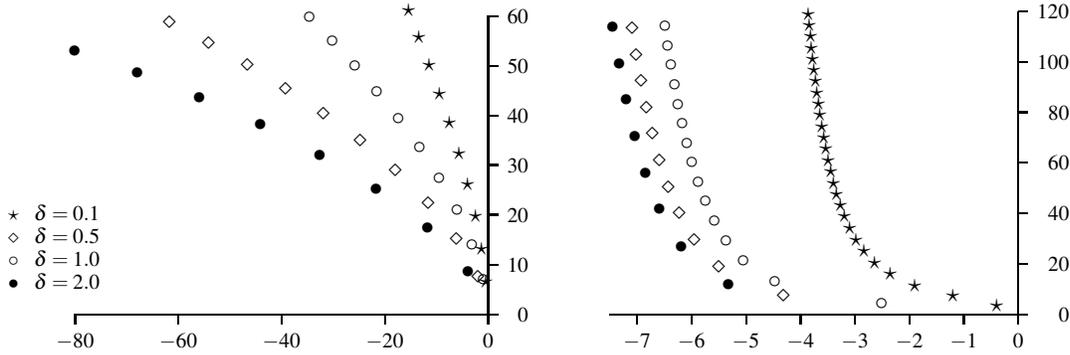

\hbox to \hsize{\hss\copy\figureten\hss\hss\copy\figureeleven\hss}
\caption{\small {\bf Roots of $\mathbf{\omega_{\,\text{\sc fz}}(s)}$ for various $\delta$
values with:\quad  left $\;\mathbf{\alpha=\frac{1}{2}}\,$, right 
$\;\mathbf{\alpha=\frac{9}{10}}$ }}
\label{fig:zener_roots}
\end{figure}

Also here we have employed the method for numerically computing roots as described above. 
In particular we use this in order to illustrate the influence of $\delta>0$ on the behaviour of the roots, see Figure~\oldref{fig:zener_roots}.

The location of the poles can thus be computed from knowledge of the constants
$b_1$, $b_2$, $c$, the exponents $\alpha$ and $\beta$ as well as the values
of $\lambda_n$.
These poles are also obtainable from the time trace measurements as the zeros
of the relaxation function which is the Laplace transform of this data $g$.
Thus assuming the spectrum $\{\lambda_n\}_{n=1}^\infty$ of $\mathcal{A}$ was
known it is perfectly reasonable that a least-squares fit could be made
to determine the damping constants contained in the term $D$ appearing in
equation~\eqref{eqn:Westervelt_init_D} and/or the wave speed $c$.
While an ill-conditioned problem, it would be particularly feasible if the
time trace data were measured at several points along an arc rather than at
a single point.
It is further conceivable that spectral information on the eigenvalues
of $\mathcal{A}$ could be determined, in particular those of the low frequency.
This in turn might be used to obtain knowledge on either a coefficient
in $\mathcal{A}$ or on the domain $\Omega$ itself as there is geometrical
information contained in the lowest few eigenvalues.
See, for example, \cite{GrebenkovNguyen:2013}.

\section{Reconstructions of $\mathbf{\kappa}$}\label{sec:reconstructions}

The boundary conditions for our test cases were homogeneous Dirichlet
at $x=0$ and homogeneous Neumann at $x=1$ with a nonhomogeneous driving term
$r(x,t)$ with greater weight near $x=1$.
Thus the solution was small in the region near $x=0$
in comparison to near $x=1$ where the data $h(t) = u(1,t)$ was measured.
The consequence of this was that $\kappa(x)$ for $x$ small was multiplied
by terms that were also small in comparison to that at the rightmost endpoint
and resulting in greater ill-conditioning of the inversion near $x=0$.
This will be apparent in each of the reconstructions to be shown below.

The data $h(t)$ was computed by the direct solver at the endpoint $x=1$
and a sample of $50$ points were taken to which uniformly distributed
random noise was added as representing the actual data measurements
that formed the overposed data.
This was then filtered by a smoothing routine based on penalising
the $H^2$ norm with regularisation parameter based on the estimate of
the noise and then up-resolved to the working size for the inverse solver.

Discretisation of $\kappa$ was done by means of a fixed set of $40$ chapeau basis functions and 
we applied a regularised frozen Newton iteration, stopped by the discrepancy principle, for numerically solving the discretised inverse problem.

\begin{figure}[h]
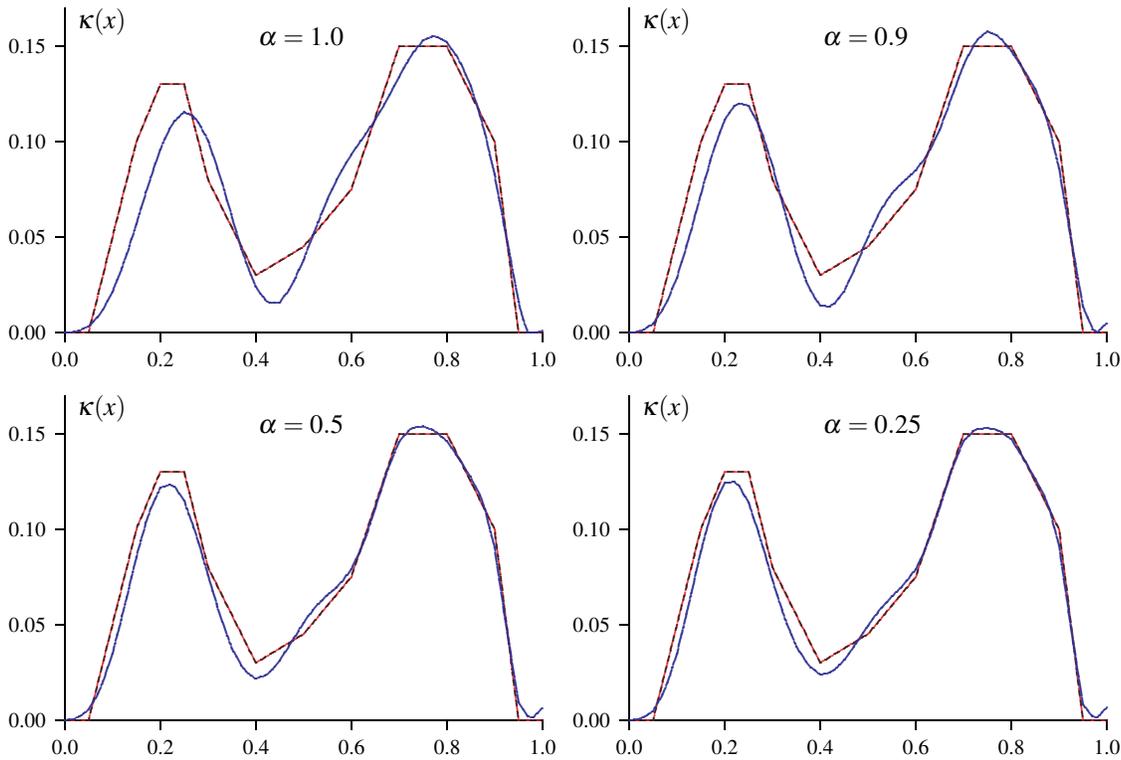

\hbox to \hsize{\hss\copy\figureone\hss\copy\figuretwo\hss}
\vskip10pt
\hbox to \hsize{\hss\copy\figurethree\hss\copy\figurefour\hss}
\caption{\small {\bf Reconstructions of $\kappa(x)$ for various $\alpha$ values.
\quad Noise = 0.1\%}}
\label{fig:kappa_with_alpha}
\end{figure}

Figure~\oldref{fig:kappa_with_alpha} shows the reconstruction of a piecewise linear $\kappa$
for the values $\alpha=1$, $\alpha=0.9$, $\alpha=0.5$ and $\alpha=0.25$.
In each case the damping coefficient  $b$ was kept at $b=0.1$ and the wave
speed at $c=1$.
The $(L^\infty,\ L^2)$ norm difference for the final versus the actual
reconstruction were: 
$(0.109, 0.078)$,
$(0.116, 0.084)$,
$(0.184, 0.126)$,
$(0.315, 0.191)$,
respectively and show the increase in resolution possible with a decrease
in $\alpha$.

Note that the reconstructions of $\kappa$ are clearly superior at the
right-hand endpoint due to imposed conditions as the wave is essentially
transmitting information primarily from right to left but the amplitude is
damped as it travels.
The smaller the fractional damping the lesser is this effect which is
also apparent from these figures.

\begin{figure}[ht]
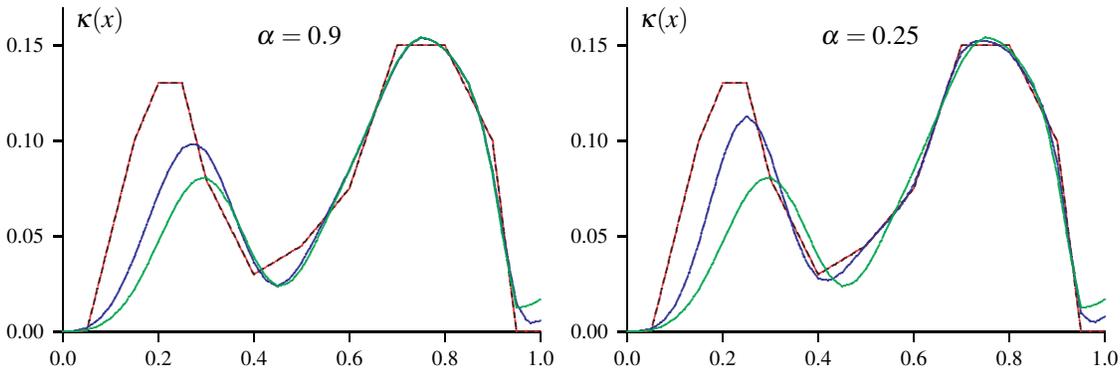

\hbox to \hsize{\hss\copy\figurethirteen\hss\copy\figuretwelve\hss}
\caption{\small {\bf Reconstructions of $\kappa(x)$ for $\alpha =0.25,\;0.9$
\quad Noise = 0.5\% (blue) and 1\% (green)}}
\label{fig:kappa_with_alpha_noise5}
\end{figure}

The reconstructions naturally worsen with increasing noise levels
as Figure~\oldref{fig:kappa_with_alpha_noise5} shows.

\begin{figure}[ht]
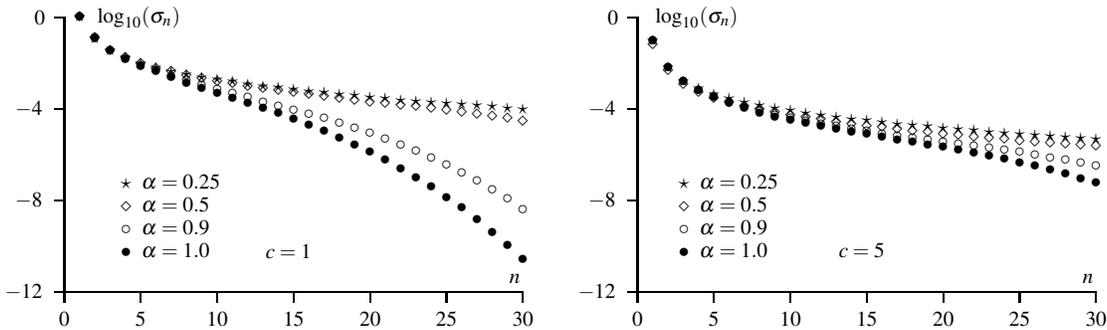

\hbox to \hsize{\hss\copy\figurefive\hss\hss\copy\figurefourteen\hss}
\caption{\small {\bf Singular values for various $\alpha$ values: left; $c=1$, right; $c=5$}}
\label{fig:sv_with_alpha}
\end{figure}

Figure~\oldref{fig:sv_with_alpha} shows the singular values of the
Jacobian matrix used in the (frozen) Newton method.
Note that if the function $\kappa$ can be well represented by a small number
of basis functions then the dependence with respect to $\alpha$ will be
fairly weak.
On the other hand, if a large number of basis functions are needed for
$\kappa$ to be represented, then the dependence on $\alpha$ becomes much
stronger although by this point the condition number of the Jacobian is
already extremely high for all $\alpha$ and relatively few singular values
are likely to be usable in any
reconstruction with data subject to extremely small noise levels.
The effect of damping is to directly contribute to the
ill-conditioning and thus it is clear that for fixed $\alpha$ and $c$ this
will increase as the coefficient $b$ increases.
The degree of ill-conditioning as a function of the wave speed $c$ is less
clear.

Figure~\oldref{fig:sv_with_alpha} shows the singular values
$\{\sigma_n\}$ of the Jacobian matrix for both $c=1$ and $c=5$.
This illustrates the decay of the singular values and hence the level of
ill-conditioning does depend on $c$ but certainly not uniformly for all
values of the fractional exponent $\alpha$.
For $\alpha$ near unity, that is damping approaches or is at
the classical paradigm, there is a considerable increase in the smaller,
high index singular values indicating the problem is much less ill-posed.
for larger wave speeds $c$.
For the smaller index $\sigma_n$ the ratio $\sigma_n/\sigma_1$ is
almost the same indicating at most a weak effect due to the wave speed.
Thus for a function $\kappa(x)$ requiring only a small number of basis functions
the effect of wave speed is relatively minimal
but this changes quite dramatically if a larger number of singular values
are required.
For $\alpha$ less than about one half the condition number $\sigma_n/\sigma_1$
becomes relatively independent of $c$ -- at least in the range indicated.

\vskip-3pt
\subsection*{Acknowledgment}
\vskip-4pt
\noindent
The work of the first author was supported by the Austrian Science Fund {\sc fwf}
under the grants P30054 and {\sc doc}78.

\noindent
The work of the second author was supported
in part by the
National Science Foundation through award {\sc dms}-1620138.

\end{document}